\documentclass{amsart}
\usepackage{amsfonts}
\usepackage{bbm}
\usepackage{mathrsfs}
\usepackage{color}
\usepackage{amsmath}
\usepackage{bigdelim}
\usepackage{multirow}
\usepackage{amssymb}
\usepackage{indentfirst}
\usepackage{CJK}
\usepackage{fancyhdr}
\usepackage{amscd,amssymb,amsmath,graphicx,verbatim,xy}
\usepackage[TS1,OT1,T1]{fontenc}
\usepackage{float}
\usepackage{booktabs}
\usepackage{mathtools}
\usepackage{bm}
\usepackage{enumitem}
\usepackage{url}

\usepackage{float}
\usepackage{booktabs}
\usepackage{mathtools}
\usepackage{bm}
\usepackage{enumitem}
\usepackage{url}
\usepackage[colorlinks=true,linkcolor=black,citecolor=black,urlcolor=black]{hyperref}

\newtheorem{theorem}{Theorem}[section]
\newtheorem{lemma}[theorem]{Lemma}
\newtheorem{corollary}[theorem]{Corollary}
\newtheorem{proposition}[theorem]{Proposition}

\theoremstyle{definition}
\newtheorem{definition}[theorem]{Definition}

\theoremstyle{remark}
\newtheorem{remark}[theorem]{Remark}
\numberwithin{equation}{section}

\makeatletter
\@namedef{subjclassname@2020}{\textup{2020} Mathematics Subject Classification}
\makeatother

\newcommand{\C}{\mathbb C}
\newcommand{\D}{\mathbb D}
\newcommand{\Dual}{\mathbb D_{\mathbb C}}
\newcommand{\B}{\overline{\mathbb B}_2}
\newcommand{\Bopen}{\mathbb B_2^{\circ}}
\newcommand{\Tset}{\mathcal T}
\newcommand{\Ptwo}{P_2(\mathfrak I)}
\newcommand{\odotop}{\mathbin{\odot}}
\newcommand{\norm}[1]{\left\lVert #1\right\rVert}
\newcommand{\abs}[1]{\left\lvert #1\right\rvert}
\newcommand{\clco}{\overline{\operatorname{conv}}}

\begin{document}
	
	\title[Multiplicatively Closed Sets in Complex Dual Numbers]{Maximal Multiplicatively Closed Subsets of the Euclidean Unit Ball in the Algebra of Complex Dual Numbers}
	
	\author{Cailing Yao}
	\address{Cailing Yao, School of Medical Informatics\\
		Changchun University of Chinese Medicine, Changchun 130117, P. R. China}
	\email{yaocl@ccucm.edu.cn}
	
	\author{Bingzhe Hou}
	\address{Bingzhe Hou, School of Mathematics, Jilin University, Changchun 130012, P. R. China}
	\email{houbz@jlu.edu.cn}
	
	\author{Chao Xie}
	\address{Chao Xie, School of Medical Informatics\\
		Changchun University of Chinese Medicine, Changchun 130117, P. R. China}
	\email{2512181948@qq.com}
	
	\author{Yuanjun Zou}
	\address{Yuanjun Zou, School of Medical Informatics\\
		Changchun University of Chinese Medicine, Changchun 130117, P. R. China}
	\email{zouyj@ccucm.edu.cn}
	
	\date{}

	\begin{abstract}
		
		We study multiplicatively closed subsets of the closed Euclidean
		unit ball in the algebra of complex dual numbers.
		An explicit maximal multiplicatively closed extension of the
		Toeplitz set \(\Tset\) is constructed.
		Every maximal multiplicatively closed subset contains the Toeplitz
		set and is compact and convex with nonempty interior.
		Such maximal subsets need not be unique.
		A finite criterion for all finite products of prescribed dual
		numbers to remain in the ball is established and applied to
		discrete-time systems under arbitrary switching, with the identity
		as initial value.
		Connections with the first two Taylor coefficients of inner
		functions and with Reinhardt geometry are also discussed.
		For a published insulin-delay subsystem, we determine the largest
		uniform input bound ensuring robust positive invariance of a
		normalized error region obtained from the constructed maximal set.
	\end{abstract}

	\maketitle
	\enlargethispage{2pt}

\section{Introduction and preliminaries}\label{sec:intro}

The concept of dual numbers dates back to W. K. Clifford's work on
biquaternions \cite{Clifford1871}. For a dual number, the coefficient of the nilpotent element
\(\varepsilon\) represents first-order information. For a polynomial, and more generally for an
analytic function \(f\), the dual evaluation
\(f(x+h\varepsilon)=f(x)+hf'(x)\varepsilon\) retains both the value and
the directional derivative. This property underlies forward-mode algorithmic
differentiation \cite{GriewankWalther2008}.
Dual-number methods also have applications in kinematics,
statics, and dynamics \cite{Fischer1999}. In control theory, unit dual quaternions provide a singularity-free and computationally compact representation of coupled attitude and position dynamics; they have been used to design sliding-mode controllers for underactuated spacecraft relative motion \cite{wu2014coupled} and neural-network-based tracking controllers for autonomous aerial vehicles \cite{zamora2025dual}. In medical imaging, dual quaternions have been applied to rigid-body motion correction in MRI, where they yield a one-step closed-form solution to the hand-eye calibration problem between an external tracking device and the scanner frame \cite{zahneisen2014fast}. In the theory of analytic functions, dual numbers arise naturally in the study of harmonic and regular functions over dual-number algebras, where a dual-number-valued harmonic function admits a regular primitive constructed via a dual-number analogue of the Cauchy integral formula \cite{deakin1966functions,kim2016dual}.

Over \(\C\), the same algebra is the two-term
quotient of analytic power series and has a lower triangular Toeplitz
representation, which connects this algebra with classical Schur
interpolation \cite{ChoCurtoLee2002,Schur1917}. These analytic and
computational connections make the interaction between multiplication
and geometric norm constraints a natural problem. A multiplicatively
closed set is preserved by multiplication by any of its elements.
This property connects the present problem with the study of invariant
sets in control theory \cite{Blanchini1999}.
In normed-algebra theory, compatibility of a norm with
multiplication is expressed by submultiplicativity
\cite{Dales2000}.
A norm \(p\) is submultiplicative if
\(p(ab)\leq p(a)p(b)\) for all algebra elements \(a,b\).
Equivalently, its closed unit ball is multiplicatively closed. For instance, the closed unit ball of \(\C\)
with the usual norm is multiplicatively closed, while that of the dual numbers is not.
This equivalence suggests studying multiplicatively closed
subsets of a prescribed unit ball when the norm itself is
not submultiplicative.

The complex dual-number algebra
\[
\Dual=\C[\varepsilon]/(\varepsilon^2)
\]
is a basic nonreduced finite-dimensional commutative complex algebra.
Each element has a unique representation \(z_0+z_1\varepsilon\),
where \(z_0,z_1\in\C\) and \(\varepsilon^2=0\). For
\(z_0+z_1\varepsilon,w_0+w_1\varepsilon\in\Dual\), multiplication is given by
\[
(z_0+z_1\varepsilon)(w_0+w_1\varepsilon)
=z_0w_0+(z_0w_1+z_1w_0)\varepsilon.
\]
For \(z=(z_0,z_1), w=(w_0,w_1)\in\C^2\),
 define the product ``$\odotop$'' in coordinates as
\begin{equation}\label{eq:product}
	z\odotop w=(z_0w_0,z_0w_1+z_1w_0).
\end{equation}
The operation ``$\odotop$'' is commutative and associative, and its identity is
\(\mathbf1=(1,0)\). We identify \(z_0+z_1\varepsilon\) with the pair
\(z=(z_0,z_1)\in\C^2\), then the complex dual-number algebra \(\Dual\) is isomorphic to \(\C^2\) with ``$\odotop$''. For convenience, we also write the original element $z_0+z_1\varepsilon$ in $\Dual$ in the form of $(z_0,z_1)$. From this perspective, the operation ``$\odotop$'' denotes the multiplication in $\Dual$. Accordingly, \((\Dual,\odotop)\) is a commutative monoid
\cite{Howie1995}. We use
\(\mathbb N=\{1,2,\ldots\}\) and
\(\mathbb N_0=\{0,1,2,\ldots\}\), and put
\(z^{\odotop0}=\mathbf1\). Direct calculation gives
\begin{equation}\label{eq:powers}
	(a,b)^{\odotop n}=(a^n,na^{n-1}b),
	\qquad n\geq1.
\end{equation}

We equip \(\Dual\cong\C^2\) with the Euclidean norm
\[
\norm{z_0+z_1\varepsilon}_2
=(\abs{z_0}^2+\abs{z_1}^2)^{1/2}.
\]
For vectors in \(\C^2\), we also write \(\norm z\) for
\(\norm z_2\). The closed unit ball is
\[
\B=\{z=(z_0,z_1)\in\C^2:
\norm z_2^2=\abs{z_0}^2+\abs{z_1}^2\leq1\},
\]
and its interior is the open unit ball
\[
\Bopen=\{z=(z_0,z_1)\in\C^2:
\abs{z_0}^2+\abs{z_1}^2<1\}.
\]
The closed unit disk \(\overline\D=\{z\in\C:\abs z\leq1\}\) is
multiplicatively closed, since \(\abs{zw}=\abs z\,\abs w\leq1\)
for \(z,w\in\overline\D\). The corresponding statement fails for the
Euclidean unit ball in \(\Dual\). For example, if
\(z=2^{-1/2}(1+\varepsilon)\), then \(\norm z_2=1\), whereas
\[
z^2=\frac12+\varepsilon,
\qquad \norm{z^2}_2^2=\frac54>1.
\]
Thus the Euclidean norm is not an algebra norm, and \(\B\) is not
multiplicatively closed. We therefore study multiplicatively closed
subsets of \(\B\), with particular attention to maximal ones.

\begin{definition}\label{def:invariant}
	For subsets \(Z,W\subseteq\Dual\), define
	\[
	Z\odotop W=\{z\odotop w:z\in Z,\ w\in W\}.
	\]
	A set \(V\subseteq\Dual\) is called \emph{multiplicatively closed} if
	\(V\odotop V\subseteq V\). If it also contains \(\mathbf1\), then it is
	a submonoid of \(\Dual\).
\end{definition}
The empty set is regarded as multiplicatively closed, and a
multiplicatively closed set is not required to contain \(\mathbf1\).

For a family \(\mathcal F\) of multiplicatively closed subsets of
\(\B\), a set \(M\in\mathcal F\) is called
\emph{maximal in \(\mathcal F\)} if there is no \(N\in\mathcal F\)
with \(M\subsetneq N\). It is called
\emph{greatest in \(\mathcal F\)} if every \(N\in\mathcal F\)
satisfies \(N\subseteq M\). Thus a maximal set has no proper extension
in the specified family, whereas a greatest set contains every set in
that family. All uses of the words \emph{maximal} and \emph{greatest}
below refer to the family stated in the surrounding theorem.

The faithful linear representation
\begin{equation}\label{eq:Toeplitz-map}
	T(z)=\begin{pmatrix}z_0&0\\z_1&z_0\end{pmatrix}
\end{equation}
satisfies
\begin{equation}\label{eq:homomorphism}
	T(z\odotop w)=T(z)T(w).
\end{equation}
We shall use the Toeplitz set
\[
\Tset=\{(z_0,z_1)\in\B:
\abs{z_1}\leq1-\abs{z_0}^2\}.
\]
Equivalently,
\[
\Tset=\{z\in\B:\norm{T(z)}_{\rm op}\leq1\},
\]
where \(\norm{\cdot}_{\rm op}\) is the operator norm induced by the
Euclidean norm on \(\C^2\). This equivalence and the basic properties
of \(\Tset\) are proved in Proposition~\ref{prop:Toeplitz-set}.

Classical coefficient interpolation provides an analytic
interpretation of the Toeplitz set \(\Tset\).
The Schur class consists of analytic functions \(f\) on the
open unit disk \(\D\) satisfying
\(\abs{f(\zeta)}\leq1\) for every \(\zeta\in\D\).
The Carath\'eodory--Fej\'er interpolation criterion
characterizes the finite sequences that occur as the initial
Taylor coefficients of Schur-class functions by the
contractivity of the associated lower triangular Toeplitz
matrices \cite{ChoCurtoLee2002,Schur1917}.
C. Yao, B. Hou and Y. Cao \cite{YaoHouCao2025} give the corresponding
characterizations for inner functions and finite Blaschke
products of bounded degree.
In particular, the first two Taylor coefficients of inner
functions form exactly \(\Tset\).

These interpolation results concern the analytic realization
of coefficient vectors.
In the dual-number algebra, \(\Tset\) is also the closed unit
ball of the algebra norm
\(z\mapsto\norm{T(z)}_{\rm op}\).
The present problem is to determine how this set can be
enlarged within the fixed Euclidean ball \(\B\) while
preserving multiplicative closure.

The remainder of the paper is organized as follows.
In Section~\ref{sec:toeplitz}, we establish the basic properties of
\(\Tset\), derive product estimates, and construct a multiplicatively
closed subset of \(\B\) that properly contains \(\Tset\).
In Section~\ref{sec:maximal}, we construct an explicit maximal
multiplicatively closed subset of \(\B\), prove the existence of
maximal extensions, and study their geometric properties and
nonuniqueness. In Section~\ref{sec:applications}, we discuss Taylor coefficients
of inner functions and Reinhardt geometry. We also give a finite
criterion for all finite products of prescribed dual numbers to
remain in \(\B\), and determine the exact input bound for an
invariant error region of a published insulin-delay subsystem.

\section{The Toeplitz set and its multiplicatively closed extensions}
\label{sec:toeplitz}

In this section, we establish the basic properties of the Toeplitz
set \(\Tset\) and construct a multiplicatively closed subset of
\(\B\) that properly contains \(\Tset\).

\subsection{Basic properties of the Toeplitz set}
\label{sec:toeplitz-properties}

The compactness and convexity of the Toeplitz contraction
set are included in \cite[Theorem~1]{ChoCurtoLee2002}.
We record the two-coordinate description and verify the
properties needed below.
\begin{proposition}\label{prop:Toeplitz-set}
	The set
	\begin{align}
		\Tset
		&:=\{z=(z_0,z_1)\in\B:
		\abs{z_1}\leq1-\abs{z_0}^2\}
		\label{eq:T-explicit}\\
		&=\{z\in\B:\norm{T(z)}_{\rm op}\leq1\}
		\label{eq:T-Toeplitz}
	\end{align}
	is compact, convex, and multiplicatively closed.
\end{proposition}
\begin{proof}
	First, the inequality
	\(\abs{z_1}\leq1-\abs{z_0}^2\) implies \(\abs{z_0}\leq1\).
	Setting \(r=\abs{z_0}\), we obtain
	\[
	\abs{z_0}^2+\abs{z_1}^2
	\leq r^2+(1-r^2)^2
	=1-r^2(1-r^2)\leq1.
	\]
	Thus the coefficient inequality already ensures that \(z\in\B\).	
	
	Let \(z=(z_0,z_1),w=(w_0,w_1)\in\Tset\), and put \(a=\abs{z_0}\),
	\(b=\abs{w_0}\).  Then
	\begin{align*}
		\abs{z_0w_1+z_1w_0}
		&\leq a(1-b^2)+b(1-a^2)\\
		&\leq1-a^2b^2,
	\end{align*}
	since the difference in the second inequality is
	\((1-a)(1-b)(1-ab)\geq0\).  Thus \(z\odotop w\in\Tset\).
	The set \(\Tset\) is closed by the continuity of
	\(z\mapsto\abs{z_0}^2+\abs{z_1}\), and it is bounded since
	\(\Tset\subseteq\B\). Hence \(\Tset\) is compact. For \(0\leq\lambda\leq1\), direct expansion gives
	\[
	\lambda\abs{z_0}^2+(1-\lambda)\abs{w_0}^2
	-\abs{\lambda z_0+(1-\lambda)w_0}^2
	=\lambda(1-\lambda)\abs{z_0-w_0}^2\geq0.
	\]
	Since \(z,w\in\Tset\), the triangle inequality therefore yields
	\begin{align*}
		\abs{\lambda z_1+(1-\lambda)w_1}
		&\leq\lambda\abs{z_1}+(1-\lambda)\abs{w_1}\\
		&\leq\lambda(1-\abs{z_0}^2)
		+(1-\lambda)(1-\abs{w_0}^2)\\
		&=1-\lambda\abs{z_0}^2-(1-\lambda)\abs{w_0}^2\\
		&\leq1-\abs{\lambda z_0+(1-\lambda)w_0}^2.
	\end{align*}
	Hence \(\lambda z+(1-\lambda)w\in\Tset\), which proves the convexity.
	
	We next prove the equivalence of the coefficient inequality
	\eqref{eq:T-explicit} and the Toeplitz contraction condition
	\eqref{eq:T-Toeplitz}. Direct calculation gives
	\[
	I-T(z)^*T(z)=
	\begin{pmatrix}
		1-\abs{z_0}^2-\abs{z_1}^2&-\overline{z_1}z_0\\
		-\overline{z_0}z_1&1-\abs{z_0}^2
	\end{pmatrix}
	\]
	and
	\[
	\det(I-T(z)^*T(z))
	=(1-\abs{z_0}^2)^2-\abs{z_1}^2.
	\]
	Since \(z\in\B\), we have
	\[
	1-\abs{z_0}^2-\abs{z_1}^2\geq0,
	\qquad 1-\abs{z_0}^2\geq0.
	\]
	Thus both diagonal entries of \(I-T(z)^*T(z)\) are nonnegative.
	The criterion for positive semidefiniteness of a \(2\times2\)
	Hermitian matrix therefore gives
	\begin{align*}
		\norm{T(z)}_{\rm op}\leq1
		&\Longleftrightarrow I-T(z)^*T(z)\text{ is positive semidefinite}\\
		&\Longleftrightarrow (1-\abs{z_0}^2)^2-\abs{z_1}^2\geq0\\
		&\Longleftrightarrow \abs{z_1}\leq1-\abs{z_0}^2.
	\end{align*}
	The last equivalence uses \(1-\abs{z_0}^2\geq0\).
	This proves the equivalence of
	\eqref{eq:T-explicit} and \eqref{eq:T-Toeplitz}.
\end{proof}
The Toeplitz set is exactly the set of dual numbers whose
multiplication preserves the whole Euclidean unit ball:
\[
\Tset
=\{w\in\Dual:\B\odotop\{w\}\subseteq\B\}.
\]
Indeed, \(z\odotop w=T(w)z\), so the condition on the
right is equivalent to \(\norm{T(w)}_{\rm op}\leq1\).

Every maximal multiplicatively closed subset of \(\B\)
contains \(\mathbf1\). Indeed, let \(V\) be such a maximal
subset. Since \(V\cup\{\mathbf1\}\) is multiplicatively
closed and contained in \(\B\), maximality gives
\(\mathbf1\in V\).
Thus the maximal sets considered below are submonoids,
although a multiplicatively closed set need not contain
\(\mathbf1\).

The following elementary facts will be used repeatedly.  In
part~(iii), 
\[
[0,1]W=\{tw:0\leq t\leq1,\ w\in W\}.
\]
\begin{proposition}\label{prop:permanence}
	Let \(Z,W\subseteq\Dual\) be multiplicatively closed.  Then
	\begin{enumerate}[label=\textup{(\roman*)}]
		\item \(Z\odotop W\) is multiplicatively closed;
		\item \(\overline W\), \(\operatorname{conv}W\), and \(\clco W\) are
		multiplicatively closed;
		\item \([0,1]W\) is multiplicatively closed.
	\end{enumerate}
\end{proposition}

\begin{proof}
	For \(z_1,z_2\in Z\) and \(w_1,w_2\in W\), associativity and commutativity
	give
	\[
	(z_1\odotop w_1)\odotop(z_2\odotop w_2)
	=(z_1\odotop z_2)\odotop(w_1\odotop w_2)\in Z\odotop W.
	\]
	This proves (i).  Since \(\odotop\) is continuous, \(\overline W\) is
	multiplicatively closed.  If
	\(u=\sum_i\lambda_iw_i\) and \(v=\sum_j\mu_jw'_j\) are finite convex
	combinations of points of \(W\), bilinearity yields
	\[
	u\odotop v
	=\sum_{i,j}\lambda_i\mu_j(w_i\odotop w'_j)
	\in\operatorname{conv}W.
	\]
	Hence \(\operatorname{conv}W\) is multiplicatively closed, and the
	closure result gives the assertion for \(\clco W\).  Finally,
	\[
	(su)\odotop(tv)=st(u\odotop v)\in[0,1]W
	\]
	for \(s,t\in[0,1]\) and \(u,v\in W\), proving (iii).
\end{proof}

\subsection{Radial coordinates and an extension of the Toeplitz set}

In this subsection, we express the ball and Toeplitz constraints in
terms of \(x=-\log\abs{z_0}\) and \(u=z_1/z_0\) for \(z_0\neq0\).
We use the resulting product estimates to construct a
multiplicatively closed subset of \(\B\) that properly contains
\(\Tset\). Throughout the paper, \(\log\) denotes the natural
logarithm.

An element \(z_0+z_1\varepsilon\in \Dual\) is invertible precisely when
\(z_0\neq0\).  In that case it factors as
\[
z_0+z_1\varepsilon=z_0(1+u\varepsilon),
\qquad u=\frac{z_1}{z_0}.
\]
For \(z=(z_0,z_1)\in\C^2\) with
\(0<\abs{z_0}\leq1\), we introduce the variables
\begin{equation}\label{eq:coordinates}
	x=-\log\abs{z_0},\qquad u=\frac{z_1}{z_0}.
\end{equation}
If \(z\) and \(w\) have parameters \((x,u)\) and \((y,v)\),
respectively, then their dual-number product \(z\odotop w\)
has parameters \((x+y,u+v)\). The variables \(x\) and \(u\) do not determine the argument of \(z_0\).
However, the defining inequalities for \(\B\) and \(\Tset\) depend
only on \(\abs{z_0}\) and \(\abs u\), so this argument is not needed
for these descriptions. Define
\begin{equation}\label{eq:GH}
	G(x)=\sqrt{e^{2x}-1},\qquad
	H(x)=e^x-e^{-x}=2\sinh x,\qquad x\geq0.
\end{equation}
For a fixed \(x\geq0\), the fiber over \(x\) consists of all
complex values of \(u=z_1/z_0\) allowed by the defining inequality
of the set. Since \(\abs{z_0}=e^{-x}\), this describes the possible
ratios \(z_1/z_0\) when the modulus of \(z_0\) is fixed.
For \(\B\) and \(\Tset\), the corresponding fibers are closed disks
in the complex \(u\)-plane, centered at the origin and with radii
\(G(x)\) and \(H(x)\), respectively.
Equivalently,
\begin{equation}\label{eq:fibers}
	z\in\B\Longleftrightarrow\abs u\leq G(x),
	\qquad
	z\in\Tset\Longleftrightarrow\abs u\leq H(x).
\end{equation}

The following estimate will be used in the extension arguments.

\begin{lemma}\label{lem:compatibility}
	For all $x,y\geq0$,
	\begin{equation}\label{eq:compatibility}
		G(x+y)\geq G(x)+H(y).
	\end{equation}
	Consequently,
	\begin{equation}\label{eq:BP}
		\B\odotop\Tset\subseteq\B.
	\end{equation}
\end{lemma}

\begin{proof}
	It follows from the definition of $G(\cdot)$ that all terms in \eqref{eq:compatibility} are nonnegative.  A direct
	calculation gives
	\begin{equation}\label{eq:square-identity}
		G(x+y)^2-(G(x)+H(y))^2
		=(e^{2y}-1)\bigl(G(x)-e^{-y}\bigr)^2\geq0.
	\end{equation}
	If \(z_0\neq0\) and \(w_0\neq0\), then
	$\abs{u+v}\leq G(x)+H(y)\leq G(x+y)$, which proves
	\eqref{eq:BP}.  If \(z_0=0\) or \(w_0=0\), the same conclusion
	follows directly from \eqref{eq:product}.
\end{proof}

\begin{lemma}\label{lem:threshold}
	If $x,y\geq\frac12\log3$, then
	\begin{equation}\label{eq:threshold}
		G(x+y)\geq G(x)+G(y).
	\end{equation}
\end{lemma}

\begin{proof}
	Expanding the difference of squares gives
	\begin{equation}\label{eq:threshold-identity}
		G(x+y)^2-(G(x)+G(y))^2
		=G(x)G(y)\bigl(G(x)G(y)-2\bigr).
	\end{equation}
	The assumption $x,y\geq\frac12\log3$ implies $G(x)\geq\sqrt2$ and $G(y)\geq\sqrt2$, so the right-hand side of \eqref{eq:threshold-identity} is
	nonnegative. Combining this with the nonnegativity of $G(\cdot)$, we obtain
	\[
	G(x+y)\geq G(x)+G(y).
	\]
\end{proof}

In the following proposition, we construct a multiplicatively closed
subset of \(\B\) that properly contains \(\Tset\).

\begin{proposition}\label{prop:old-extension}
	The multiplicatively closed set \(\Tset\) is not maximal in \(\B\).
	More precisely, there exists \(\delta_0>0\) such that,
	for every real number \(\delta\) satisfying
	\(0<\delta<\min\{\delta_0,e^{-2}-e^{-4}\}\),
	the following construction gives a strict extension.  Put
	\begin{equation}\label{eq:xi}
		\xi=\left(e^{-1},1-e^{-2}+\delta\right),
		\qquad
		\Lambda_\xi=\{\mathbf1,\xi,\xi^{\odotop2},\ldots\}.
	\end{equation}
	Then
	\begin{equation}\label{eq:old-extension}
		\Tset\subsetneq \Lambda_\xi\odotop\Tset\subseteq\B.
	\end{equation}
\end{proposition}

\begin{proof}
	Put \(\widetilde\xi=(e^{-1},1-e^{-2})\). By
	\eqref{eq:coordinates}, its parameters are \(x=1\) and
	\(u=e-e^{-1}=H(1)<G(1)\).
	Let \(w=(w_0,w_1)\in\Tset\) with \(w_0\neq0\), and set
	\[
	y=-\log\abs{w_0},\qquad v=\frac{w_1}{w_0}.
	\]
	It follows from Lemma~\ref{lem:compatibility} that
	\[
	\abs{u+v}\leq H(1)+H(y)<G(1)+H(y)\leq G(1+y).
	\]
	If \(w_0=0\), then \(w=(0,w_1)\) with \(\abs{w_1}\leq1\). By the multiplication formula,
	\[
	\widetilde\xi\odotop w=(0,e^{-1}w_1),
	\]
	and hence
	\[
	\norm{\widetilde\xi\odotop w}
	=e^{-1}\abs{w_1}\leq e^{-1}<1.
	\]
	Therefore, \(\widetilde\xi\odotop w\in\Bopen\) for every
	\(w=(0,w_1)\in\Tset\).
	
	Since \(\Tset\) is compact and
	\(\widetilde\xi\odotop w\in\Bopen\) for every \(w\in\Tset\),
	\[
	m:=\max_{w\in\Tset}\norm{\widetilde\xi\odotop w}_2<1.
	\]
	Set \(\delta_0=1-m>0\).
	For \(\abs{\delta}<\delta_0\) and \(w=(w_0,w_1)\in\Tset\),
	the multiplication formula yields
	\[
	\left(e^{-1},1-e^{-2}+\delta\right)\odotop w
	=\widetilde\xi\odotop w+(0,\delta w_0).
	\]
	Since \(\abs{w_0}\leq1\), the triangle inequality shows that
	\[
	\norm{\left(e^{-1},1-e^{-2}+\delta\right)\odotop w}
	\leq \norm{\widetilde\xi\odotop w}
	+\abs{\delta}\abs{w_0}
	\leq m+\abs{\delta}<1.
	\]
	Therefore,
	\[
	\left\{\left(e^{-1},1-e^{-2}+\delta\right)\right\}
	\odotop\Tset\subseteq\Bopen,
	\quad\text{whenever }\abs{\delta}<\delta_0.
	\]
	Choose a real number $\delta$ such that
	\[
	0<\delta<\min\{\delta_0,e^{-2}-e^{-4}\}.
	\]
	Then
	\begin{equation}\label{eq:xi-T-open}
		\{\xi\}\odotop\Tset\subseteq\Bopen,
		\qquad \xi\notin\Tset.
	\end{equation}
	
	By \eqref{eq:powers}, for $n\geq2$,
	\[
	\xi^{\odotop n}
	=\left(e^{-n},ne^{-(n-1)}(1-e^{-2}+\delta)\right).
	\]
	Since $ne^{-(n-1)}\leq1$ for $n\geq2$ and
	$1-e^{-2}+\delta<1-e^{-4}\leq1-e^{-2n}$, we have
	\begin{equation}\label{eq:xi-higher-powers}
		\xi^{\odotop n}\in\Tset,
		\quad n\geq2.
	\end{equation}
	The set \(\Lambda_\xi\) is multiplicatively
	closed, so it follows from Proposition~\ref{prop:permanence} that
	\(\Lambda_\xi\odotop\Tset\) is multiplicatively closed.  
	Since \(\Tset\) is multiplicatively closed,
	\eqref{eq:xi-higher-powers} gives
	\[
	\xi^{\odotop n}\odotop\Tset\subseteq\Tset,
	\quad n\geq2.
	\]
	Together with \eqref{eq:xi-T-open} and
	\(\xi^{\odotop0}=\mathbf1\), this implies
	\(\Lambda_\xi\odotop\Tset\subseteq\B\).
	Since \(\mathbf1\in\Lambda_\xi\), we have
	\(\Tset\subseteq\Lambda_\xi\odotop\Tset\).
	This inclusion is strict because
	\(\xi=\xi\odotop\mathbf1\in\Lambda_\xi\odotop\Tset\),
	whereas \(\xi\notin\Tset\). Therefore, \(\Tset\) is not maximal among the multiplicatively
	closed subsets of \(\B\).
\end{proof}

\section{Maximal multiplicatively closed subsets}\label{sec:maximal}

In this section, we construct an explicit maximal multiplicatively
closed subset of \(\B\), prove the existence of maximal extensions,
and establish their geometric properties. We also show that there
exist at least two distinct maximal multiplicatively closed subsets
of \(\B\), and that no multiplicatively closed subset of \(\B\)
contains every other such subset.

\subsection{A multiplicatively closed extension of the Toeplitz set}

In Proposition~\ref{prop:old-extension}, we constructed a
multiplicatively closed extension of \(\Tset\) by adjoining
\(\xi\) and taking all finite products of elements of
\(\Tset\cup\{\xi\}\).
In this subsection, we construct a multiplicatively closed set
\(E\) containing \(\Lambda_\xi\odotop\Tset\), which will be used
to obtain an explicit maximal multiplicatively closed subset
of \(\B\).
Let
\begin{equation}\label{eq:E}
	L=\{z\in\B:\abs{z_0}\leq3^{-1/2}\},
	\qquad E=\Tset\cup L.
\end{equation}

\begin{theorem}\label{thm:explicit-extension}
	The set \(E\) is multiplicatively closed and
	\begin{equation}\label{eq:E-strict}
		\Tset\subsetneq E\subseteq\B.
	\end{equation}
	Moreover, \(\Lambda_\xi\odotop\Tset\subseteq E\).
\end{theorem}

\begin{proof}
	Since \(\Tset\) is multiplicatively closed, we have
	\[
	\Tset\odotop\Tset\subseteq\Tset\subseteq E.
	\]
	Let \(z=(z_0,z_1)\in L\) and
	\(w=(w_0,w_1)\in\Tset\). By \eqref{eq:BP},
	\(z\odotop w\in\B\). Moreover,
	\[
	\abs{z_0w_0}
	\leq\abs{z_0}
	\leq3^{-1/2},
	\]
	since \(\abs{w_0}\leq1\). Thus \(z\odotop w\in L\), and
	commutativity yields
	\[
	L\odotop\Tset=\Tset\odotop L\subseteq L\subseteq E.
	\]
	
	Now let \(z=(z_0,z_1)\in L\) and \(w=(w_0,w_1)\in L\),
	with \(z_0\neq0\) and \(w_0\neq0\). Since \(z,w\in L\), we have
	\(\abs{z_0},\abs{w_0}\leq3^{-1/2}\). Hence
	\[
	x=-\log\abs{z_0}\geq\tfrac12\log3,
	\qquad
	y=-\log\abs{w_0}\geq\tfrac12\log3.
	\]
	Set \(u=z_1/z_0\) and \(v=w_1/w_0\). It follows from Lemma~\ref{lem:threshold} that
	\[
	\abs{u+v}\leq G(x)+G(y)\leq G(x+y).
	\]
	Therefore,
	\begin{align*}
		\norm{z\odotop w}_2^2
		&=e^{-2(x+y)}\bigl(1+\abs{u+v}^2\bigr)\\
		&\leq e^{-2(x+y)}\bigl(1+G(x+y)^2\bigr)=1,
	\end{align*}
	which shows that \(z\odotop w\in\B\). The first coordinate of \(z\odotop w\) is \(z_0w_0\), and
	\[
	\abs{z_0w_0}\leq3^{-1/2}\cdot3^{-1/2}=\frac13
	\leq3^{-1/2}.
	\]
	Together with \(z\odotop w\in\B\), this gives
	\(z\odotop w\in L\).  If \(z_0=0\), then
	\[
	z\odotop w=(0,z_1w_0),
	\qquad
	\norm{z\odotop w}_2
	=\abs{z_1w_0}\leq3^{-1/2}<1,
	\]
	since \(\abs{z_1}\leq1\) and \(\abs{w_0}\leq3^{-1/2}\).
	
	If \(w_0=0\), then
	\[
	z\odotop w=(0,z_0w_1),
	\qquad
	\norm{z\odotop w}_2
	=\abs{z_0w_1}\leq3^{-1/2}<1.
	\]
	Thus \(z\odotop w\in\B\) and its first coordinate is zero
	whenever \(z_0=0\) or \(w_0=0\). By the definition of \(L\),
	we obtain \(z\odotop w\in L\).  Thus \(E\) is
	multiplicatively closed.
	
	To prove that \(\Tset\subsetneq E\), consider
	\begin{equation}\label{eq:q-first}
		q=\left(3^{-1/2},\sqrt{2/3}\right).
	\end{equation}
	Since
	\[
	\norm q_2^2=\frac13+\frac23=1,
	\qquad \abs{q_0}=3^{-1/2},
	\]
	we have \(q\in L\subseteq E\). However,
	\[
	\abs{q_1}=\sqrt{2/3}>\frac23=1-\abs{q_0}^2,
	\]
	so \(q\notin\Tset\). Therefore, \(\Tset\subsetneq E\).
	
	Finally, since \(\mathbf1\in\Tset\),
	\eqref{eq:xi-T-open} implies \(\xi\in\B\). Together with
	\[
	\abs{\xi_0}=e^{-1}<3^{-1/2},
	\]
	this gives \(\xi\in L\subseteq E\).
	Since \(E\) is multiplicatively closed and contains \(\Tset\),
	we obtain
	\[
	\xi^{\odotop n}\odotop w\in E
	\qquad(n\in\mathbb N_0,\ w\in\Tset).
	\]
	Consequently,
	\[
	\Lambda_\xi\odotop\Tset\subseteq E.
	\]
\end{proof}

\begin{remark}\label{rem:threshold-sharp}
	For \(0<r\leq1\), let
	\[
	L_r=\{z=(z_0,z_1)\in\B:\abs{z_0}\leq r\}.
	\]
	Then \(L_r\) is multiplicatively closed if and only if
	\(r\leq3^{-1/2}\).
	
	Suppose first that \(r\leq3^{-1/2}\). Then \(L_r\subseteq L\).
	The proof of Theorem~\ref{thm:explicit-extension} shows that
	\(L\) is multiplicatively closed. Thus, for
	\(z=(z_0,z_1),w=(w_0,w_1)\in L_r\), we have
	\(z\odotop w\in L\subseteq\B\). Moreover,
	\[
	\abs{z_0w_0}\leq r^2\leq r.
	\]
	Hence \(z\odotop w\in L_r\), proving that \(L_r\) is
	multiplicatively closed.
	
	Conversely, suppose that \(r>3^{-1/2}\). Choose \(s\) such that
	\[
	3^{-1/2}<s<1,\qquad s\leq r,
	\]
	and put \(z=(s,\sqrt{1-s^2})\). Since
	\[
	\norm z_2^2=s^2+(1-s^2)=1,
	\qquad \abs{z_0}=s\leq r,
	\]
	we have \(z\in L_r\). However,
	\[
	z\odotop z=(s^2,2s\sqrt{1-s^2}),
	\]
	and therefore
	\[
	\norm{z\odotop z}_2^2-1
	=4s^2-3s^4-1
	=(1-s^2)(3s^2-1)>0.
	\]
	Thus \(z\odotop z\notin\B\), so \(L_r\) is not
	multiplicatively closed.
	
	Consequently, \(3^{-1/2}\) is the largest parameter \(r\)
	for which \(L_r\) is multiplicatively closed.
	The greatest multiplicatively closed subset of \(\B\)
	containing \(E\) will be determined in the next subsection.
\end{remark}

\subsection{An explicit maximal multiplicatively closed subset}

In this subsection, we describe subsets of \(\B\) by prescribing
a bound on \(\abs{z_1/z_0}\) for each fixed value of
\(x=-\log\abs{z_0}\), where \(z_0\neq0\).
We then determine the greatest multiplicatively closed subset
of \(\B\) containing \(E\).

Let \(f:[0,\infty)\to[0,\infty)\) be any function satisfying
\[
0\leq f(x)\leq G(x)=\sqrt{e^{2x}-1},
\quad x\geq0.
\]
Define 
\begin{align}
	\mathcal R_f
	:=\;&\{(z_0,z_1)\in\B:z_0\neq0,\ 
	\abs{z_1/z_0}\leq f(-\log\abs{z_0})\}\notag\\
	&{}\cup\{(0,z_1):\abs{z_1}\leq1\}.
	\label{eq:radial-set}
\end{align}
For each fixed \(x\geq0\), the ratios \(u=z_1/z_0\)
corresponding to points of \(\mathcal R_f\) with
\(\abs{z_0}=e^{-x}\) form the closed disk
\[
\{u\in\C:\abs u\leq f(x)\},
\]
centered at the origin. For \(z_0=0\), the definition includes
all points \((0,z_1)\) with \(\abs{z_1}\leq1\).
These points form the fiber of \(\mathcal R_f\) over \(z_0=0\).

We say that \(f\) is \emph{superadditive} if
\begin{equation}\label{eq:superadditive}
	f(x+y)\geq f(x)+f(y),
	\quad x,y\geq0.
\end{equation}
\begin{proposition}\label{prop:radial-invariance}
	The set \(\mathcal R_f\) is multiplicatively closed if and only
	if \(f\) is superadditive.
\end{proposition}

\begin{proof}
	First, we suppose that $f$ is superadditive. Let \(z=(z_0,z_1)\in\mathcal R_f\) and
	\(w=(w_0,w_1)\in\mathcal R_f\), with
	\(z_0\neq0\) and \(w_0\neq0\). Set
	\[
	x=-\log\abs{z_0},\qquad u=\frac{z_1}{z_0},
	\qquad
	y=-\log\abs{w_0},\qquad v=\frac{w_1}{w_0}.
	\]
	By \eqref{eq:product}, the parameters of \(z\odotop w\) satisfy
	\[
	-\log\abs{z_0w_0}=x+y,
	\qquad
	\frac{z_0w_1+z_1w_0}{z_0w_0}=u+v.
	\] By \eqref{eq:superadditive} and the triangle
	inequality,
	\[
	\abs{u+v}\leq f(x)+f(y)\leq f(x+y).
	\]
	Since \(f\leq G\), it follows that
	\[
	\norm{z\odotop w}_2^2
	=e^{-2(x+y)}\bigl(1+\abs{u+v}^2\bigr)
	\leq e^{-2(x+y)}\bigl(1+G(x+y)^2\bigr)=1.
	\]
	Thus \(z\odotop w\in\mathcal R_f\).
	If \(z_0=0\), then \eqref{eq:product} shows that
	\[
	z\odotop w=(0,z_1w_0),
	\qquad
	\abs{z_1w_0}=\abs{z_1}\abs{w_0}\leq1,
	\]
	since \(z,w\in\mathcal R_f\subseteq\B\).
	If \(w_0=0\), then
	\[
	z\odotop w=(0,z_0w_1),
	\qquad
	\abs{z_0w_1}=\abs{z_0}\abs{w_1}\leq1.
	\]
	By the definition of \(\mathcal R_f\), we have
	\(z\odotop w\in\mathcal R_f\) whenever \(z_0=0\) or \(w_0=0\).
	
	We have shown that, if \(f\) is superadditive, then
	\(z\odotop w\in\mathcal R_f\) for every
	\(z,w\in\mathcal R_f\). Hence \(\mathcal R_f\) is
	multiplicatively closed.
	
	Conversely, suppose that \(\mathcal R_f\) is multiplicatively
	closed. For arbitrary \(x,y\geq0\), put
	\[
	z=\bigl(e^{-x},e^{-x}f(x)\bigr),
	\qquad
	w=\bigl(e^{-y},e^{-y}f(y)\bigr).
	\]
	Since \(0\leq f\leq G\), we have
	\begin{align*}
		\norm z_2^2
		&=e^{-2x}\bigl(1+f(x)^2\bigr)
		\leq e^{-2x}\bigl(1+G(x)^2\bigr)=1,\\
		\norm w_2^2
		&=e^{-2y}\bigl(1+f(y)^2\bigr)
		\leq e^{-2y}\bigl(1+G(y)^2\bigr)=1.
	\end{align*}
	Moreover,
	\[
	\frac{z_1}{z_0}=f(x),
	\qquad
	\frac{w_1}{w_0}=f(y).
	\]
	Thus \(z,w\in\mathcal R_f\). By \eqref{eq:product},
	\[
	z\odotop w
	=\bigl(e^{-(x+y)},
	e^{-(x+y)}(f(x)+f(y))\bigr).
	\]
	Multiplicative closure implies \(z\odotop w\in\mathcal R_f\). The first coordinate of \(z\odotop w\) is
	\(e^{-(x+y)}>0\). Applying the inequality
	\[
	\abs{p_1/p_0}\leq f(-\log\abs{p_0})
	\]
	from \eqref{eq:radial-set} to
	\[
	p=z\odotop w
	=\bigl(e^{-(x+y)},e^{-(x+y)}(f(x)+f(y))\bigr),
	\]
	we obtain
	\[
	f(x)+f(y)
	=\abs{\frac{e^{-(x+y)}(f(x)+f(y))}{e^{-(x+y)}}}
	\leq f\bigl(-\log e^{-(x+y)}\bigr)
	=f(x+y).
	\]
	Therefore, \(f\) is superadditive.
\end{proof}

We now express \(E\) in terms of the variable \(x\) and define
a function that will determine its greatest multiplicatively
closed extension in \(\B\). Let
\begin{equation}\label{eq:t0}
	t_0=\frac12\log3,\qquad G(t_0)=\sqrt2.
\end{equation}
The set \(E\) in \eqref{eq:E} is \(\mathcal R_{f_E}\), where
\[
f_E(x)=
\begin{cases}
	H(x),&0\leq x<t_0,\\
	G(x),&x\geq t_0.
\end{cases}
\]
Define the larger profile
\begin{equation}\label{eq:fmax}
	f_E^{\max}(x)=
	\begin{cases}
		G(x+t_0)-G(t_0),&0\leq x\leq t_0,\\
		G(x),&x\geq t_0.
	\end{cases}
\end{equation}
The two branches agree at \(t_0\), because
\(G(2t_0)=2G(t_0)\).

In the following lemma, we prove that \(f_E^{\max}\) is
superadditive.

\begin{lemma}\label{lem:fmax-superadditive}
	For all \(x,y\geq0\),
	\[
	f_E^{\max}(x+y)\geq f_E^{\max}(x)+f_E^{\max}(y).
	\]
\end{lemma}
\begin{proof}
	Recall that \(t_0=\frac12\log3\), and define
	\[
	J(t)=G(t+t_0)-G(t_0),
	\qquad 0\leq t\leq t_0.
	\]
	Then \(J(0)=0\) and \(J(t_0)=G(t_0)\), since
	\(G(2t_0)=2G(t_0)\). For \(t\geq t_0\), we have
	\(e^{2t}\geq3\), and therefore
	\[
	G'(t)=\frac{e^{2t}}{\sqrt{e^{2t}-1}},
	\qquad
	G''(t)=
	\frac{e^{2t}(e^{2t}-2)}{(e^{2t}-1)^{3/2}}
	\geq0.
	\]
	Thus \(G'\) is nondecreasing on \([t_0,\infty)\).
	Consequently, \(J'\) is nondecreasing and \(J\) is convex
	on \([0,t_0]\).
	
	By symmetry in \(x\) and \(y\), it suffices to consider
	\(0\leq x\leq y\). We distinguish four cases.
	
	\begin{enumerate}[label=\textup{(\roman*)}]
		\item \(0\leq x\leq y<t_0\) and \(x+y\leq t_0\).
		
		In this case, \(f_E^{\max}\) agrees with \(J\) at
		\(x\), \(y\), and \(x+y\). Since \(J'\) is nondecreasing,
		\begin{align*}
			f_E^{\max}(x+y)-f_E^{\max}(x)
			&=J(x+y)-J(x)\\
			&=\int_0^y J'(x+t)\,dt\\
			&\geq\int_0^y J'(t)\,dt\\
			&=J(y)=f_E^{\max}(y).
		\end{align*}
		
		\item \(0\leq x\leq y<t_0\) and \(x+y>t_0\).
		
		Put \(s=x+y\), so \(t_0<s<2t_0\). Define
		\[
		\phi(t)=J(t)+J(s-t),
		\qquad s-t_0\leq t\leq t_0.
		\]
		The function \(\phi\) is convex, and
		\[
		\phi(s-t_0)=\phi(t_0)
		=J(s-t_0)+J(t_0)=G(s).
		\]
		Hence \(\phi(t)\leq G(s)\) for every
		\(t\in[s-t_0,t_0]\). Since \(x\in[s-t_0,t_0]\),
		we obtain
		\[
		f_E^{\max}(x)+f_E^{\max}(y)
		=J(x)+J(y)
		=\phi(x)
		\leq G(s)
		=f_E^{\max}(x+y).
		\]
		
		\item \(0\leq x<t_0\leq y\).
		
		In this case, we have \(f_E^{\max}(x)=J(x)\),
		\(f_E^{\max}(y)=G(y)\) and
		\(f_E^{\max}(x+y)=G(x+y)\).
		Since \(G'\) is nondecreasing on \([t_0,\infty)\), we obtain
		\begin{align*}
			f_E^{\max}(x+y)-f_E^{\max}(y)
			&=G(x+y)-G(y)\\
			&=\int_0^x G'(y+t)\,dt\\
			&\geq\int_0^x G'(t_0+t)\,dt\\
			&=J(x)=f_E^{\max}(x).
		\end{align*}
		
		\item \(t_0\leq x\leq y\).
		
		In this case, \(f_E^{\max}\) agrees with \(G\) at
		\(x\), \(y\), and \(x+y\). By  Lemma~\ref{lem:threshold},
		\[
		f_E^{\max}(x+y)
		=G(x+y)
		\geq G(x)+G(y)
		=f_E^{\max}(x)+f_E^{\max}(y).
		\]
	\end{enumerate}
	
	Thus the inequality holds for all \(x,y\geq0\), and
	\(f_E^{\max}\) is superadditive.
\end{proof}

In the following lemma, we compare \(f_E^{\max}\) with
\(f_E\) and \(G\).

\begin{lemma}\label{lem:fmax-bounds}
	For every \(x\geq0\),
	\[
	f_E(x)\leq f_E^{\max}(x)\leq G(x).
	\]
	Moreover,
	\[
	f_E(x)<f_E^{\max}(x),
	\quad 0<x<t_0.
	\]
\end{lemma}

\begin{proof}
	For \(x\geq t_0\), the definitions give
	\[
	f_E(x)=f_E^{\max}(x)=G(x).
	\]
	It remains to consider \(0\leq x<t_0\), where
	\[
	f_E(x)=H(x),
	\qquad
	f_E^{\max}(x)=G(x+t_0)-G(t_0).
	\]
	
	Substituting \(t_0\) and \(x\) for the two variables
	in \eqref{eq:square-identity}, we obtain
	\[
	G(x+t_0)^2-\bigl(G(t_0)+H(x)\bigr)^2
	=(e^{2x}-1)\bigl(\sqrt2-e^{-x}\bigr)^2
	\geq0.
	\]
	Since \(G(x+t_0)\geq0\) and \(G(t_0)+H(x)\geq0\), we obtain
	\[
	G(x+t_0)\geq G(t_0)+H(x).
	\]
	Consequently,
	\[
	f_E^{\max}(x)
	=G(x+t_0)-G(t_0)
	\geq H(x)=f_E(x).
	\]
	For \(0<x<t_0\), we have \(e^{2x}-1>0\) and
	\(\sqrt2-e^{-x}>0\). Therefore,
	\[
	f_E^{\max}(x)
	=G(x+t_0)-G(t_0)
	>H(x)=f_E(x).
	\]
	
	To prove the upper bound, observe that
	\(0\leq G(x)\leq G(t_0)=\sqrt2\).
	Substituting \(y=t_0\) into
	\eqref{eq:threshold-identity} yields
	\[
	G(x+t_0)^2-\bigl(G(x)+G(t_0)\bigr)^2
	=\sqrt2\,G(x)\bigl(\sqrt2\,G(x)-2\bigr)
	\leq0.
	\]
	Therefore,
	\[
	G(x+t_0)\leq G(x)+G(t_0),
	\]
	and hence
	\[
	f_E^{\max}(x)
	=G(x+t_0)-G(t_0)
	\leq G(x).
	\]
	This completes the proof.
\end{proof}

Set
\[
M_E^{\max}=\mathcal R_{f_E^{\max}},
\]
and define
\begin{equation}\label{eq:rho-max}
	\rho_{\max}(r)=
	\begin{cases}
		\sqrt{1-r^2},
		&0\leq r\leq3^{-1/2},\\
		\sqrt{3-r^2}-\sqrt2\,r,
		&3^{-1/2}\leq r\leq1.
	\end{cases}
\end{equation}
For \(0<r\leq1\), substitution of \(x=-\log r\) gives
\[
r f_E^{\max}(-\log r)=\rho_{\max}(r).
\]
For \(z_0=0\), the definition of \(\mathcal R_{f_E^{\max}}\)
includes exactly the points \((0,z_1)\) with
\(\abs{z_1}\leq1=\rho_{\max}(0)\). Hence
\begin{equation}\label{eq:explicit-maximal}
	M_E^{\max}
	=\{(z_0,z_1)\in\C^2:
	\abs{z_0}\leq1,\ 
	\abs{z_1}\leq\rho_{\max}(\abs{z_0})\}.
\end{equation}
At \(r=3^{-1/2}\), both expressions for \(\rho_{\max}\)
equal \(\sqrt{2/3}\). Moreover,
\(\rho_{\max}(0)=1\) and \(\rho_{\max}(1)=0\).

In the following theorem, we prove that \(M_E^{\max}\)
is the greatest multiplicatively closed subset of \(\B\)
containing \(E\).

\begin{theorem}\label{thm:explicit-maximal}
	The set \(M_E^{\max}\) is the greatest multiplicatively
	closed subset of \(\B\) containing \(E\).
	It is closed and convex, and satisfies
	\[
	\Tset\subsetneq E\subsetneq M_E^{\max}\subseteq\B.
	\]
	In particular, \(M_E^{\max}\) is maximal among all
	multiplicatively closed subsets of \(\B\).
\end{theorem}

\begin{proof}
	By Lemmas~\ref{lem:fmax-superadditive}
	and~\ref{lem:fmax-bounds}, the function \(f_E^{\max}\) is
	superadditive and satisfies \(0\leq f_E^{\max}\leq G\). Therefore, it follows from
	Proposition~\ref{prop:radial-invariance} that
	\(M_E^{\max}\) is multiplicatively closed.
	The inequality \(f_E\leq f_E^{\max}\) gives
	\(E\subseteq M_E^{\max}\). For every \(x\in(0,t_0)\),
	\[
	\bigl(e^{-x},e^{-x}f_E^{\max}(x)\bigr)
	\in M_E^{\max}\setminus E,
	\]
	since \(f_E(x)<f_E^{\max}(x)\).
	Together with Theorem~\ref{thm:explicit-extension}, this yields
	\[
	\Tset\subsetneq E\subsetneq M_E^{\max}\subseteq\B.
	\]
	
	Let \(a=3^{-1/2}\), and write the two expressions for
	\(\rho_{\max}\) as
	\[
	\begin{aligned}
		\rho_1(r)&=\sqrt{1-r^2},&&0\leq r\leq a,\\
		\rho_2(r)&=\sqrt{3-r^2}-\sqrt2\,r,&&a\leq r\leq1.
	\end{aligned}
	\]
	Since \(\rho_1(a)=\rho_2(a)=\sqrt{2/3}\),
	the function \(\rho_{\max}\) is continuous on \([0,1]\).
	For \(0<r<a\) and \(a<r<1\), respectively,
	\[
	\begin{array}{ll}
		\rho_1'(r)=-\dfrac{r}{\sqrt{1-r^2}},&
		\rho_1''(r)=-\dfrac1{(1-r^2)^{3/2}},\\[7pt]
		\rho_2'(r)=-\dfrac{r}{\sqrt{3-r^2}}-\sqrt2,&
		\rho_2''(r)=-\dfrac3{(3-r^2)^{3/2}}.
	\end{array}
	\]
	Thus \(\rho_{\max}'(r)<0\) on
	\((0,a)\cup(a,1)\), so continuity at \(a\) implies that
	\(\rho_{\max}\) is nonincreasing on \([0,1]\).
	Moreover,
	\[
	\rho_{\max,-}'(a)=-\frac1{\sqrt2}
	>
	-\frac5{2\sqrt2}=\rho_{\max,+}'(a).
	\]
	The negative second derivatives imply that the derivative is
	nonincreasing on each of \((0,a)\) and \((a,1)\).
	For \(0<r<a<s<1\), we also have
	\[
	\rho_{\max}'(r)
	\geq\rho_{\max,-}'(a)
	>
	\rho_{\max,+}'(a)
	\geq\rho_{\max}'(s).
	\]
	Hence \(\rho_{\max}\) is concave on \([0,1]\).
	The unequal one-sided derivatives show that it is not
	differentiable at \(a\).
	
	For later use, let \(\rho:[0,1]\to[0,\infty)\) satisfy
	\(\rho(r)\leq\sqrt{1-r^2}\), and write
	\begin{equation}\label{eq:radial-profile-set}
		K_\rho=
		\{(z_0,z_1)\in\C^2:
		\abs{z_0}\leq1,\ 
		\abs{z_1}\leq\rho(\abs{z_0})\}.
	\end{equation}
	By \eqref{eq:explicit-maximal} and the convexity
	criterion proved in Lemma~\ref{lem:Reinhardt-convexity}
	below, \(M_E^{\max}\) is convex.
	Since \(\rho_{\max}\) is continuous on \([0,1]\), the conditions
	\[
	\abs{z_0}\leq1,\qquad
	\abs{z_1}\leq\rho_{\max}(\abs{z_0})
	\]
	are preserved under limits in \(\C^2\).
	Hence \(M_E^{\max}\) is closed in \(\C^2\).
	
	We now prove that every multiplicatively closed set
	\(V\subseteq\B\) containing \(E\) satisfies
	\(V\subseteq M_E^{\max}\).
	Fix such a set \(V\), and take \(z=(z_0,z_1)\in V\).
	When \(z_0\neq0\), set
	\[
	x=-\log\abs{z_0}\geq0,\qquad u=\frac{z_1}{z_0}.
	\]
	We consider the following three cases.
	
	\begin{enumerate}[label=\textup{(\roman*)}]
		\item \(z_0=0\).
		
		Since \(z\in V\subseteq\B\),
		\[
		\abs{z_1}\leq1=\rho_{\max}(0).
		\]
		Therefore \(z=(0,z_1)\in M_E^{\max}\).
		
		\item \(z_0\neq0\) and \(x\geq t_0\).
		
		The inequality \(\abs{z_0}^2+\abs{z_1}^2\leq1\) shows
		\[
		\abs u^2
		=\frac{\abs{z_1}^2}{\abs{z_0}^2}
		\leq\frac{1-\abs{z_0}^2}{\abs{z_0}^2}
		=e^{2x}-1=G(x)^2.
		\]
		Hence \(\abs u\leq G(x)=f_E^{\max}(x)\), and therefore
		\(z\in M_E^{\max}\).
		
		\item \(z_0\neq0\) and \(0\leq x<t_0\).
		
		If \(u=0\), then
		\[
		\abs u=0
		\leq G(x+t_0)-G(t_0)
		=f_E^{\max}(x),
		\]
		because \(G\) is nondecreasing. Thus \(z\in M_E^{\max}\).
		
		Now suppose \(u\neq0\). Set
		\[
		\eta=\frac{u}{\abs u},\qquad
		w=\bigl(e^{-t_0},e^{-t_0}G(t_0)\eta\bigr).
		\]
		Then \(\abs\eta=1\), \(u=\abs u\,\eta\), and
		\[
		\abs{w_0}=e^{-t_0}=3^{-1/2},\qquad
		\norm w_2^2=e^{-2t_0}\bigl(1+G(t_0)^2\bigr)=1.
		\]
		Therefore \(w\in L\subseteq E\subseteq V\).
		Since \(V\) is multiplicatively closed,
		\(z\odotop w\in V\subseteq\B\).
		By \eqref{eq:product},
		\[
		\frac{z_0w_1+z_1w_0}{z_0w_0}
		=u+G(t_0)\eta,\qquad
		\abs{z_0w_0}=e^{-(x+t_0)}>0.
		\]
		The inequality
		\(\abs{z_0w_0}^2+\abs{z_0w_1+z_1w_0}^2\leq1\)
		therefore gives
		\[
		\abs{u+G(t_0)\eta}^2
		\leq\abs{z_0w_0}^{-2}-1
		=e^{2(x+t_0)}-1
		=G(x+t_0)^2.
		\]
		Since \(u=\abs u\,\eta\) and \(\abs\eta=1\),
		\[
		\abs u+G(t_0)
		=\abs{u+G(t_0)\eta}
		\leq G(x+t_0).
		\]
		Consequently,
		\[
		\abs u
		\leq G(x+t_0)-G(t_0)
		=f_E^{\max}(x),
		\]
		so \(z\in M_E^{\max}\).
	\end{enumerate}
	
	Thus \(z\in M_E^{\max}\) for every \(z\in V\), and
	\(V\subseteq M_E^{\max}\).
	Therefore \(M_E^{\max}\) is the greatest multiplicatively
	closed subset of \(\B\) containing \(E\).
	
	Finally, let \(N\subseteq\B\) be multiplicatively closed
	with \(M_E^{\max}\subseteq N\). Then
	\[
	E\subseteq M_E^{\max}\subseteq N\subseteq\B.
	\]
	Taking \(V=N\) in the inclusion just proved yields
	\(N\subseteq M_E^{\max}\). Hence \(N=M_E^{\max}\),
	which proves that \(M_E^{\max}\) is maximal among all
	multiplicatively closed subsets of \(\B\).
\end{proof}

\begin{figure}[H]
	\centering
	\includegraphics[width=0.4\linewidth]{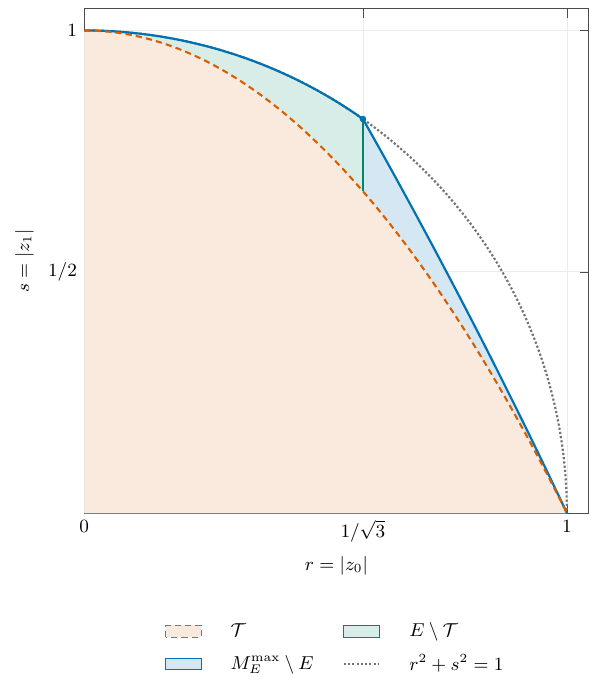}
	\caption{The sets \(\Tset\subsetneq E\subsetneq M_E^{\max}\)
		in modulus coordinates. The orange, green and blue regions
		represent \(\Tset\), \(E\setminus\Tset\) and
		\(M_E^{\max}\setminus E\), respectively.
		The vertical segment at \(r=1/\sqrt3\) is part of the
		boundary of the region representing \(E\).
		The dotted arc is \(r^2+s^2=1\).}
	\label{fig:explicit-maximal-regions}
\end{figure}

Figure~\ref{fig:explicit-maximal-regions} illustrates the inclusions
\(\Tset\subsetneq E\subsetneq M_E^{\max}\).
Each set is represented in the coordinates
\((r,s)=(\abs{z_0},\abs{z_1})\).

\subsection{Existence and geometry of maximal multiplicatively closed subsets}
In this subsection, we prove that every multiplicatively closed
subset of \(\B\) is contained in a maximal multiplicatively
closed subset of \(\B\).
We then show that every maximal multiplicatively closed subset
of \(\B\) contains \(\Tset\) and is compact and convex with
nonempty interior. Finally, we determine the Hausdorff dimensions
of these sets and their boundaries.

In the following theorem, we establish the existence of maximal
multiplicatively closed extensions in \(\B\).

\begin{theorem}\label{thm:zorn}
	Let \(A\subseteq\B\) be multiplicatively closed.
	Then there exists a maximal multiplicatively closed subset
	\(M_A\) of \(\B\) such that \(A\subseteq M_A\).
\end{theorem}

\begin{proof}
	Consider
	\[
	\mathcal F_A=
	\{V\subseteq\B:A\subseteq V,\ V\odotop V\subseteq V\},
	\]
	ordered by inclusion. Since \(A\in\mathcal F_A\), we have
	\(\mathcal F_A\neq\varnothing\).
	
	Let \(\mathcal C\) be a nonempty chain in \(\mathcal F_A\),
	and put
	\[
	U=\bigcup_{V\in\mathcal C}V.
	\]
	Then \(A\subseteq U\subseteq\B\).
	For \(z,w\in U\), choose \(V_1,V_2\in\mathcal C\) such that
	\(z\in V_1\) and \(w\in V_2\).
	Since \(\mathcal C\) is ordered by inclusion, we may assume
	\(V_1\subseteq V_2\). Thus
	\[
	z,w\in V_2,
	\qquad z\odotop w\in V_2\subseteq U.
	\]
	Hence \(U\odotop U\subseteq U\), so \(U\in\mathcal F_A\).
	Moreover, \(V\subseteq U\) for every \(V\in\mathcal C\).
	Therefore \(U\) is an upper bound of \(\mathcal C\)
	in \(\mathcal F_A\).
	
	By Zorn's lemma \cite{DaveyPriestley2002}, there exists
	\(M_A\in\mathcal F_A\) that is maximal under inclusion.
	Let \(N\subseteq\B\) be multiplicatively closed and satisfy
	\(M_A\subseteq N\). Then
	\[
	A\subseteq M_A\subseteq N\subseteq\B,
	\]
	so \(N\in\mathcal F_A\). Maximality of \(M_A\) in
	\(\mathcal F_A\) gives \(N=M_A\).
	Thus \(M_A\) is a maximal multiplicatively closed subset
	of \(\B\) containing \(A\).
\end{proof}

For \(A=E\), it follows from Theorem~\ref{thm:explicit-maximal} that
\(M_A\subseteq M_E^{\max}\).
Since \(M_A\) is maximal and \(M_E^{\max}\) is
multiplicatively closed, we obtain \(M_A=M_E^{\max}\).
Thus \(E\) has a unique maximal multiplicatively closed
extension in \(\B\), and
\begin{equation}\label{eq:E-maximal}
	\Tset\subsetneq E\subsetneq M_E^{\max}\subseteq\B.
\end{equation}

\begin{proposition}\label{thm:closed-convex}
	Let \(M\) be a maximal multiplicatively closed subset of \(\B\).
	Then \(M\) is closed and convex in \(\C^2\).
\end{proposition}

\begin{proof}
	Since \(\B\) is closed and convex, we have
	\[
	M\subseteq\clco M\subseteq\B.
	\]
	By Proposition~\ref{prop:permanence}\textup{(ii)},
	the set \(\clco M\) is multiplicatively closed.
	Since \(M\) is maximal among the multiplicatively closed
	subsets of \(\B\), we obtain
	\[
	M=\clco M.
	\]
	Thus \(M\) is closed and convex.
\end{proof}

In the following proposition, we prove that every maximal
multiplicatively closed subset of \(\B\) contains \(\Tset\).

\begin{proposition}\label{prop:automatic-Toeplitz}
	Let \(M\) be a maximal multiplicatively closed subset of \(\B\).
	Then
	\[
	M\odotop\Tset=M.
	\]
	In particular, \(\Tset\subseteq M\).
\end{proposition}

\begin{proof}
	Since \(M\cup\{\mathbf1\}\) is multiplicatively closed
	and contained in \(\B\), maximality of \(M\) implies
	\(\mathbf1\in M\).
	Together with \(\mathbf1\in\Tset\), this yields
	\[
	M\cup\Tset\subseteq M\odotop\Tset.
	\]
	
	By Proposition~\ref{prop:permanence}\textup{(i)},
	the set \(M\odotop\Tset\) is multiplicatively closed.
	Moreover, it follows from Lemma~\ref{lem:compatibility} that
	\[
	M\odotop\Tset
	\subseteq\B\odotop\Tset
	\subseteq\B.
	\]
	Thus \(M\odotop\Tset\) is a multiplicatively closed
	subset of \(\B\) containing \(M\).
	The maximality of \(M\) therefore gives
	\[
	M\odotop\Tset=M.
	\]
	Since \(\Tset\subseteq M\odotop\Tset\), we finally obtain
	\(\Tset\subseteq M\).
\end{proof}

A set \(K\subseteq\C^2\) is called \emph{complex balanced} if
\[
\lambda K\subseteq K,
\quad\lambda\in\C,\ \abs{\lambda}\leq1.
\]
A \emph{convex body} in \(\C^2\cong\mathbb R^4\) is a
compact convex set with nonempty interior.

We recall the definition of Hausdorff dimension
\cite{Mattila1995}. For \(A\subseteq\C^2\), \(s>0\),
and \(\delta>0\), define
\[
\mathcal H^s_\delta(A)
=
\inf_{\substack{
		A\subseteq\bigcup_{j=1}^{\infty}U_j\\
		\operatorname{diam}U_j\leq\delta\ \text{for all }j}}
\sum_{j=1}^{\infty}
\bigl(\operatorname{diam}U_j\bigr)^s,
\]
where the infimum is taken over countable covers of \(A\)
by subsets \(U_j\subseteq\C^2\), and
\(\operatorname{diam}U_j\) denotes the diameter with respect
to the Euclidean metric. The \(s\)-dimensional Hausdorff
measure of \(A\) is
\[
\mathcal H^s(A)
=\lim_{\delta\downarrow0}\mathcal H^s_\delta(A).
\]
The Hausdorff dimension of \(A\) is defined by
\[
\dim_{\rm H}A
=\inf\{s>0:\mathcal H^s(A)=0\}.
\]

We recall the following standard dimension result for boundaries
of convex bodies \cite{Mattila1995,Schneider2014}.

\begin{lemma}\label{lem:convex-boundary-dimension}
	Let \(K\) be a convex body in \(\mathbb R^n\), where \(n\geq2\).
	Then
	\[
	\dim_{\rm H}(\partial K)=n-1.
	\]
\end{lemma}
\begin{corollary}\label{cor:dimension}
	Let \(M\) be a maximal multiplicatively closed subset of \(\B\).
	Then \(M\) is a complex-balanced convex body, and
	\begin{equation}\label{eq:dimensions}
		\dim_{\rm H}M=4,
		\qquad
		\dim_{\rm H}(\partial M)=3,
	\end{equation}
	where \(\partial M\) is the boundary of \(M\) in \(\C^2\).
\end{corollary}

\begin{proof}
	It follows from Propositions~\ref{thm:closed-convex}
	and~\ref{prop:automatic-Toeplitz} that \(M\) is closed, convex and \(\Tset\subseteq M\).
	
	Let \(z=(z_0,z_1)\in\C^2\) satisfy
	\(\norm z_2^2\leq3/4\), and put \(t=\abs{z_1}\).
	Then
	\[
	\abs{z_0}^2+\abs{z_1}
	\leq\frac34-t^2+t
	=1-\left(t-\frac12\right)^2
	\leq1.
	\]
	Thus \(z\in\Tset\), and therefore
	\begin{equation}\label{eq:common-inner-ball}
		\frac{\sqrt3}{2}\B
		\subseteq\Tset\subseteq M\subseteq\B.
	\end{equation}
	Since \((\sqrt3/2)\B\subseteq M\), the set \(M\)
	has nonempty interior. Moreover, \(M\) is closed and bounded,
	so it is compact. Hence \(M\) is a convex body.
	
	For every \(\lambda\in\C\) with \(\abs{\lambda}\leq1\),
	we have \((\lambda,0)\in\Tset\subseteq M\). Consequently,
	\[
	\lambda z=(\lambda,0)\odotop z\in M,
	\quad z\in M.
	\]
	Thus \(M\) is complex balanced.
	
	By \eqref{eq:common-inner-ball} and the monotonicity of
	Hausdorff dimension,
	\[
	4=\dim_{\rm H}\!\left(\frac{\sqrt3}{2}\B\right)
	\leq\dim_{\rm H}M
	\leq\dim_{\rm H}\B=4.
	\]
	Hence \(\dim_{\rm H}M=4\).
	Since \(M\) is a convex body in \(\mathbb R^4\),
	it follows from Lemma~\ref{lem:convex-boundary-dimension} that
	\[
	\dim_{\rm H}(\partial M)=3.\qedhere
	\]
\end{proof}

For a maximal multiplicatively closed subset \(M\) of \(\B\),
define its Minkowski functional by
\[
p_M(z)=\inf\{t>0:z\in tM\},\qquad z\in\C^2.
\]

\begin{corollary}\label{cor:maximal-algebra-norm}
	The function \(p_M\) is a submultiplicative norm on \(\Dual\),
	with closed unit ball \(M\), and
	\[
	\norm z_2\leq p_M(z)\leq\norm{T(z)}_{\rm op},
	\qquad p_M(\mathbf1)=1.
	\]
	If \(q\) is a submultiplicative norm on \(\Dual\) such that
	\[
	\norm z_2\leq q(z)\leq p_M(z)
	\qquad\text{for every }z\in\Dual,
	\]
	then \(q=p_M\).
\end{corollary}

\begin{proof}
	Since \(M\) is a complex-balanced convex body with
	\(0\) in its interior, \(p_M\) is a norm with closed unit
	ball \(M\). The inclusions \(\Tset\subseteq M\subseteq\B\)
	give the two norm bounds and \(p_M(\mathbf1)=1\).
	For nonzero \(z,w\), multiplicative closure gives
	\[
	\frac{z\odotop w}{p_M(z)p_M(w)}
	=\frac{z}{p_M(z)}\odotop\frac{w}{p_M(w)}\in M.
	\]
	Thus \(p_M(z\odotop w)\leq p_M(z)p_M(w)\).
	The same inequality holds if either factor is zero.
	Finally, the closed unit ball \(B_q\) of \(q\) satisfies
	\[
	M\subseteq B_q\subseteq\B.
	\]
	It is multiplicatively closed, so maximality of \(M\)
	implies \(B_q=M\). Hence \(q=p_M\).
\end{proof}

\subsection{Powers of dual numbers and nonuniqueness of maximal subsets}
In this subsection, we give a criterion for all powers of an
element \(s\in\B\) to remain in \(\B\).
We then construct \(p,q\in\B\) whose powers all remain in
\(\B\), but satisfy \(p\odotop q\notin\B\), and use these
elements to prove nonuniqueness of maximal multiplicatively
closed subsets.

For \(s=(a,b)\in\B\), define
\begin{equation}\label{eq:cyclic-set}
	\Lambda_s=\{s^{\odotop n}:n\in\mathbb N_0\}.
\end{equation}
Since
\[
s^{\odotop0}=\mathbf1,\qquad
s^{\odotop m}\odotop s^{\odotop n}
=s^{\odotop(m+n)},
\quad m,n\in\mathbb N_0,
\]
the set \(\Lambda_s\) is a submonoid of \(\Dual\).

\begin{proposition}\label{prop:cyclic}
	Let \(s=(a,b)\in\B\).
	If \(a=0\), then \(\Lambda_s\subseteq\B\).
	If \(a\neq0\), then \(\Lambda_s\subseteq\B\) if and only if
	\begin{equation}\label{eq:cyclic-criterion}
		\abs a^{2n}\left(1+n^2\abs{b/a}^2\right)\leq1,
		\quad n\geq1.
	\end{equation}
\end{proposition}

\begin{proof}
	If \(a=0\), then \(s^{\odotop n}=(0,0)\) for every \(n\geq2\).
	Since \(s\in\B\) and \(\mathbf1\in\B\), we obtain
	\(\Lambda_s\subseteq\B\).
	
	If \(a\neq0\), it follows from \eqref{eq:powers} that
	\[
	\norm{s^{\odotop n}}_2^2
	=\abs a^{2n}\left(1+n^2\abs{b/a}^2\right),
	\quad n\geq1.
	\]
	Since \(s^{\odotop0}=\mathbf1\in\B\),
	we have \(\Lambda_s\subseteq\B\) if and only if
	\eqref{eq:cyclic-criterion} holds for every \(n\geq1\).
\end{proof}

Define \(m:[0,\infty)\to[0,\infty)\) by
\begin{equation}\label{eq:m-function}
	m(x)=\inf_{n\geq1}\frac{G(nx)}{n},
	\qquad x\geq0.
\end{equation}

\begin{corollary}\label{cor:cyclic-profile}
	Let \(s=(a,b)\in\B\) with \(a\neq0\).
	Then \(\Lambda_s\subseteq\B\) if and only if
	\[
	\abs{b/a}\leq m(-\log\abs a).
	\]
\end{corollary}

\begin{proof}
	Put \(x=-\log\abs a\geq0\).
	It follows from Proposition~\ref{prop:cyclic} that
	\(\Lambda_s\subseteq\B\) if and only if
	\[
	e^{-2nx}\left(1+n^2\abs{b/a}^2\right)\leq1,
	\quad n\geq1.
	\]
	For each \(n\geq1\), this inequality is equivalent to
	\[
	\abs{b/a}\leq
	\frac{\sqrt{e^{2nx}-1}}{n}
	=\frac{G(nx)}{n}.
	\]
	Consequently, \(\Lambda_s\subseteq\B\) if and only if
	\[
	\abs{b/a}\leq
	\inf_{n\geq1}\frac{G(nx)}{n}=m(x).
	\]
\end{proof}

For \(x>0\), the value of \(m(x)\) can be determined from
finitely many terms. To see this, set
\[
g_x(t)=\frac{G(tx)}{t},
\qquad t>0.
\]
For \(t\geq x^{-1}\), logarithmic differentiation gives
\[
\frac{d}{dt}\log g_x(t)
=\frac{x}{1-e^{-2tx}}-\frac1t
\geq x-\frac1t\geq0.
\]
Thus \(g_x\) is nondecreasing on \([x^{-1},\infty)\).
Let \(N_x=\lceil x^{-1}\rceil\), where \(\lceil x^{-1}\rceil\)
denotes the smallest integer greater than or equal to \(x^{-1}\).
Then
\[
m(x)=\min_{1\leq n\leq N_x}\frac{G(nx)}{n}.
\]
For \(x=0\), we have \(G(0)=0\), so \(m(0)=0\).

\begin{corollary}\label{cor:generated-extension}
	If $\Lambda_s\subseteq\B$, then
	\begin{equation}\label{eq:Vs}
		V_s=\Lambda_s\odotop\Tset
	\end{equation}
	is a multiplicatively closed subset of \(\B\) containing
	\(\Tset\cup\{s\}\).
\end{corollary}

\begin{proof}
	It follows from Proposition~\ref{prop:permanence} that
	\(V_s=\Lambda_s\odotop\Tset\) is multiplicatively closed.
	Since \(\Lambda_s\subseteq\B\), it can be deduced from
	Lemma~\ref{lem:compatibility} that \(V_s\subseteq\B\).
	Moreover, \(\mathbf1\in\Lambda_s\cap\Tset\), so
	\[
	\Tset=\{\mathbf1\}\odotop\Tset\subseteq V_s,
	\qquad
	s=s\odotop\mathbf1\in V_s.
	\]
	Thus \(\Tset\cup\{s\}\subseteq V_s\).
\end{proof}

We now exhibit two dual numbers whose cyclic submonoids are each
contained in \(\B\), but which cannot belong to a common multiplicatively
closed subset of \(\B\).  In coordinate form, put
\begin{equation}\label{eq:pq}
	p=\left(\frac34,\frac{21}{40}\right),
	\qquad
	q=\left(\frac1{\sqrt3},\sqrt{\frac23}\right).
\end{equation}
Both points $p$ and $q$ belong to $\B\setminus\Tset$.  Formula \eqref{eq:powers}
gives
\begin{equation}\label{eq:p-powers}
	\norm{p^{\odotop n}}_2^2
	=\left(\frac9{16}\right)^n
	\left(1+\frac{49n^2}{100}\right)=:A_n.
\end{equation}
Direct calculation gives \(A_1,A_2,A_3<1\), with
\(A_3=394389/409600\). For $n\geq3$, the inequality
$A_{n+1}<A_n$ is equivalent to
\[
343n^2-882n+259>0,
\]
which holds because the polynomial is increasing for $n\geq3$ and has
value $700$ at $n=3$.  Therefore $\Lambda_p\subseteq\B$.  Similarly,
\begin{equation}\label{eq:q-powers}
	\norm{q^{\odotop n}}_2^2=3^{-n}(1+2n^2)\leq1.
\end{equation}
The inequality \(1+2n^2\leq3^n\) holds at \(n=1\).
The induction step follows from
\[
3(1+2n^2)-\bigl(1+2(n+1)^2\bigr)
=4n(n-1)\geq0,\qquad n\geq1.
\]
Hence \(\Lambda_q\subseteq\B\).
On the other hand,
\[
p\odotop q=
\left(\frac{\sqrt3}{4},\frac{\sqrt6}{4}+\frac{7\sqrt3}{40}\right),
\]
and hence
\begin{equation}\label{eq:incompatible}
	\begin{aligned}
		\norm{p\odotop q}_2^2
		&=\frac3{16}+
		\left(\frac{\sqrt6}{4}+\frac{7\sqrt3}{40}\right)^2\\
		&=\frac{1047}{1600}+\frac{21\sqrt2}{80}>1.
	\end{aligned}
\end{equation}
Therefore \(p\odotop q\notin\B\).
Consequently, no multiplicatively closed subset of \(\B\)
contains both \(p\) and \(q\).

The following theorem establishes the nonuniqueness of maximal
multiplicatively closed subsets of \(\B\).

\begin{theorem}\label{thm:nonunique}
	There exist maximal multiplicatively closed subsets \(M_1\) and
	\(M_2\) of \(\B\) such that
	\[
	\Tset\subseteq M_1\cap M_2,
	\qquad
	M_1\nsubseteq M_2,
	\qquad
	M_2\nsubseteq M_1.
	\]
	Neither the family of all multiplicatively closed subsets of \(\B\)
	nor its subfamily consisting of sets containing \(\Tset\)
	has a greatest element.
\end{theorem}

\begin{proof}
	By Corollary~\ref{cor:generated-extension}, the set \(V_p\)
	is multiplicatively closed and satisfies
	\[
	\Tset\cup\{p\}\subseteq V_p\subseteq\B.
	\]
	It follows from Theorem~\ref{thm:zorn} that there exists a maximal
	multiplicatively closed subset \(M_1\) of \(\B\) containing \(V_p\).
	Put \(M_2=M_E^{\max}\).
	By Theorem~\ref{thm:explicit-maximal}, the set \(M_2\) is maximal.
	Since \(\Tset\cup\{q\}\subseteq E\subseteq M_2\), we have
	\[
	\Tset\subseteq M_1\cap M_2,
	\qquad p\in M_1,
	\qquad q\in M_2.
	\]
	
	Multiplicative closure gives
	\[
	\begin{aligned}
		q\in M_1
		&\Longrightarrow p\odotop q\in M_1\subseteq\B,\\
		p\in M_2
		&\Longrightarrow p\odotop q\in M_2\subseteq\B.
	\end{aligned}
	\]
	Each conclusion contradicts \eqref{eq:incompatible}.
	Therefore
	\[
	p\in M_1\setminus M_2,
	\qquad
	q\in M_2\setminus M_1.
	\]
	In particular, neither \(M_1\) nor \(M_2\) is contained in the other.
	
	Finally, assume that a greatest set \(M_*\) exists in either
	of the two families.
	Since both \(M_1\) and \(M_2\) contain \(\Tset\), the defining
	property of \(M_*\) shows that
	\[
	M_1\cup M_2\subseteq M_*\subseteq\B.
	\]
	Thus \(p,q\in M_*\), and multiplicative closure implies
	\[
	p\odotop q\in M_*\odotop M_*
	\subseteq M_*\subseteq\B,
	\]
	again contradicting \eqref{eq:incompatible}.
	Hence neither family has a greatest element.
\end{proof}

\begin{remark}\label{rem:relative-greatest}
	The set \(M_E^{\max}\) is greatest among the multiplicatively
	closed subsets of \(\B\) containing \(E\), by
	Theorem~\ref{thm:explicit-maximal}.
	The larger family of multiplicatively closed subsets of \(\B\)
	containing \(\Tset\) has no greatest element, by
	Theorem~\ref{thm:nonunique}.
\end{remark}

\section{Connections and applications}\label{sec:applications}
In this section, we discuss connections with inner functions and
Reinhardt geometry. We also apply our results to discrete-time
switched systems and an insulin-delay subsystem.

\subsection{The first two Taylor coefficients of inner functions}

In this subsection, we recall the identification of
\(\Tset\) with the set of pairs formed by the first
two Taylor coefficients of inner functions. Classical results on the Taylor coefficients of inner functions
were obtained by D. J. Newman and H. S. Shapiro
\cite{NewmanShapiro1962}.
Let \(\D\) be the open unit disk, and let \(\mathfrak I\) denote
the class of inner functions on \(\D\), including unimodular
constant functions.  Let \(f\in\mathfrak I\), and write its Taylor expansion
at the origin as
\[
f(\zeta)=c_0+c_1\zeta+c_2\zeta^2+\cdots,
\qquad \zeta\in\D.
\]
Define
\[
P_2(f)=(c_0,c_1)=(f(0),f'(0)).
\] 
Taking products modulo \(\zeta^2\) identifies the coefficient pair
\((c_0,c_1)\) with the dual number \(c_0+c_1\varepsilon\).
The representation \(T\) then connects the
dual-number problem with operator theory and analytic function theory.
The Carath\'eodory--Fej\'er interpolation theorem characterizes
the initial Taylor coefficients of Schur-class functions by
contractivity of the associated lower triangular Toeplitz matrix
\cite{ChoCurtoLee2002,Schur1917}. C.~Yao, B.~Hou and Y.~Cao 
\cite{YaoHouCao2025} give the corresponding formulations for inner
functions and for finite Blaschke products of bounded degree.  The identity \(\Ptwo=\Tset\) is the case \(n=2\) of
\cite[Theorem~2.5]{YaoHouCao2025}.
The following statement records this identification
together with the properties of \(\Tset\).

\begin{theorem}\label{thm:coefficient-body}
	The set \(\Ptwo\) has the following descriptions:
	\begin{align}
		\Ptwo
		&=\{(f(0),f'(0)):f\in\mathfrak I\}=\Tset,
		\label{eq:P-inner}\\
		&=\{z\in\B:\norm{T(z)}_{\rm op}\leq1\},
		\label{eq:P-Toeplitz}\\
		&=\{z=(z_0,z_1)\in\B:
		\abs{z_1}\leq1-\abs{z_0}^2\}.
		\label{eq:P-explicit}
	\end{align}
	Moreover, \(\Ptwo\) is compact, convex, and multiplicatively closed.
\end{theorem}

\begin{proof}
	By \cite[Theorem~2.5]{YaoHouCao2025} with \(n=2\),
	\[
	\Ptwo
	=\{z\in\C^2:\norm{T(z)}_{\rm op}\leq1\}
	=\Tset.
	\]
	The remaining assertions follow from
	Proposition~\ref{prop:Toeplitz-set}.
\end{proof}

\subsection{Reinhardt sets and radial boundary functions}

In this subsection, we describe the rotational symmetry and
convexity of the sets defined by radial boundary functions.
A set \(K\subseteq\C^2\) is called a \emph{Reinhardt set}
\cite{JarnickiPflug2008} if
\[
(z_0,z_1)\in K
\quad\Longrightarrow\quad
(e^{i\theta_0}z_0,e^{i\theta_1}z_1)\in K,
\quad \theta_0,\theta_1\in\mathbb R.
\]
The imaginary unit is denoted by \(i\).
The angles \(\theta_0\) and \(\theta_1\) may be chosen separately.
Complex balance guarantees invariance under a common rotation
\((z_0,z_1)\mapsto(e^{i\theta}z_0,e^{i\theta}z_1)\),
but does not in general imply the Reinhardt property.

The sets \(K_\rho\) defined in
\eqref{eq:radial-profile-set} are contained in \(\B\)
and are Reinhardt sets.
For each fixed \(z_0\), the allowed values of \(z_1\)
form the closed disk centered at the origin with radius
\(\rho(\abs{z_0})\).
The following elementary criterion describes convexity in terms
of \(\rho\).

\begin{lemma}\label{lem:Reinhardt-convexity}
	The set \(K_\rho\) is convex if and only if \(\rho\) is
	concave and nonincreasing on \([0,1]\).
\end{lemma}

\begin{proof}
	Suppose first that \(\rho\) is concave and nonincreasing.
	Let \(z,w\in K_\rho\) and \(0\leq\lambda\leq1\).
	Then
	\begin{align*}
		\abs{\lambda z_1+(1-\lambda)w_1}
		&\leq\lambda\rho(\abs{z_0})+(1-\lambda)\rho(\abs{w_0})\\
		&\leq\rho(\lambda\abs{z_0}+(1-\lambda)\abs{w_0})\\
		&\leq\rho(\abs{\lambda z_0+(1-\lambda)w_0}).
	\end{align*}
	Also, \(\abs{\lambda z_0+(1-\lambda)w_0}\leq1\), so
	\(\lambda z+(1-\lambda)w\in K_\rho\).
	
	Conversely, suppose that \(K_\rho\) is convex.
	For \(r,s\in[0,1]\), both \((r,\rho(r))\) and
	\((s,\rho(s))\) belong to \(K_\rho\). Hence
	\[
	\lambda\rho(r)+(1-\lambda)\rho(s)
	\leq\rho(\lambda r+(1-\lambda)s),
	\quad 0\leq\lambda\leq1,
	\]
	which proves concavity.
	For \(0\leq r\leq s\leq1\) with \(s>0\), put
	\(\lambda=(1+r/s)/2\). Since
	\((s,\rho(s)),(-s,\rho(s))\in K_\rho\), convexity gives
	\[
	(r,\rho(s))
	=\lambda(s,\rho(s))+(1-\lambda)(-s,\rho(s))\in K_\rho.
	\]
	Thus \(\rho(s)\leq\rho(r)\), proving that \(\rho\) is
	nonincreasing.
\end{proof}

The explicit set \(M_E^{\max}\) is \(K_{\rho_{\max}}\).
Its convexity is therefore expressed by the concavity and
monotonicity of \(\rho_{\max}\), as verified in the proof of
Theorem~\ref{thm:explicit-maximal}.
For \(0<r\leq1\), the change of variables
\[
x=-\log r,\qquad \rho(r)=r f(-\log r)
\]
converts the inequality \(\abs{z_1}\leq\rho(\abs{z_0})\)
into \(\abs{z_1/z_0}\leq f(x)\).
When \(\rho(0)=1\), the set \(K_\rho\) is exactly
\(\mathcal R_f\). Proposition~\ref{prop:radial-invariance}
then characterizes its multiplicative closure by superadditivity
of \(f\). Thus convexity and multiplicative closure are described
by two corresponding one-variable functions.

\subsection{Finite criteria for products of prescribed dual numbers}

In this subsection, we give a finite criterion for all finite
products of prescribed dual numbers to remain in \(\B\).
The asymptotic growth of matrix products and the associated
extremal norms are studied in joint spectral-radius theory
\cite{Wirth2002}.
Our criterion concerns the Euclidean norms of the finite
products themselves, including products of small length.
We also interpret the condition for a discrete-time system
that starts at \(\mathbf1\) and is multiplied at each step
by one of the prescribed elements.
Let \(S=\{s_1,\ldots,s_d\}\subseteq\B\), where
\(s_j=(a_j,b_j)\), and consider
\begin{equation}\label{eq:switching-system}
	\zeta_{t+1}=\zeta_t\odotop s_{\sigma_t},\qquad
	\zeta_0=\mathbf1,\qquad
	t\in\mathbb N_0,\quad \sigma_t\in\{1,\ldots,d\}.
\end{equation}
At each step, any element of \(S\) may be chosen as the next
factor, and repetitions are allowed. The sequence
\(\sigma=(\sigma_t)_{t\in\mathbb N_0}\) is the switching signal.
The system is called \emph{safe under arbitrary switching from
	the identity} if \(\zeta_t\in\B\) for every
\(t\in\mathbb N_0\) and every switching signal.
Its reachable set is
\[
\mathcal R_S=\{\mathbf1\}\cup
\{s_{i_1}\odotop\cdots\odotop s_{i_k}:
k\in\mathbb N,\ i_\ell\in\{1,\ldots,d\}\}.
\]
For each switching signal \(\sigma\), the recursion gives
\[
\zeta_t=s_{\sigma_0}\odotop\cdots\odotop s_{\sigma_{t-1}}
\qquad(t\geq1).
\]
Conversely, every finite product of elements of \(S\)
is obtained by a suitable choice of the switching signal. Thus safety is equivalent to
\(\mathcal R_S\subseteq\B\).
The following consequence of the maximal-extension theorem
expresses this condition in terms of maximal multiplicatively
closed subsets.

\begin{corollary}\label{cor:safety-maximal}
	Let \(S\subseteq\B\) be finite. Then
	\(\mathcal R_S\subseteq\B\) if and only if there is a maximal
	multiplicatively closed subset \(M\) of \(\B\) such that
	\(S\subseteq M\).
\end{corollary}

\begin{proof}
	Assume that \(\mathcal R_S\subseteq\B\).
	For \(x,y\in\mathcal R_S\), write
	\[
	x=s_{i_1}\odotop\cdots\odotop s_{i_m},
	\qquad
	y=s_{j_1}\odotop\cdots\odotop s_{j_n},
	\]
	where \(m,n\geq0\), with the empty product interpreted
	as \(\mathbf1\). Then
	\[
	x\odotop y
	=s_{i_1}\odotop\cdots\odotop s_{i_m}
	\odotop s_{j_1}\odotop\cdots\odotop s_{j_n}
	\in\mathcal R_S.
	\]
	Thus \(\mathcal R_S\) is multiplicatively closed. It follows from
	Theorem~\ref{thm:zorn} that there exists a maximal
	multiplicatively closed subset \(M\) of \(\B\) such that
	\(\mathcal R_S\subseteq M\). Since each element of \(S\) is a product with a single
	factor, we obtain \(S\subseteq\mathcal R_S\subseteq M\).
	
	Conversely, let \(M\) be a maximal multiplicatively closed
	subset of \(\B\) containing \(S\). By
	Proposition~\ref{prop:automatic-Toeplitz}, we know that 
	\(\mathbf1\in\Tset\subseteq M\). Multiplicative closure then
	implies
	\[
	s_{i_1}\odotop\cdots\odotop s_{i_k}\in M,
	\quad k\in\mathbb N,\ i_\ell\in\{1,\ldots,d\}.
	\]
	Hence \(\mathcal R_S\subseteq M\subseteq\B\).
\end{proof}

If \(S\subseteq M\), multiplicative closure also gives
\[
z\odotop s_j\in M,
\quad z\in M,\ 1\leq j\leq d.
\]
Thus trajectories starting at any point of \(M\) remain in
\(M\). This does not assert that trajectories starting at every
point of \(\B\) remain in \(\B\).

We first consider products with a zero first coordinate in one
or more factors.

\begin{lemma}\label{lem:zero-scalar}
	Let \(u_\ell=(\alpha_\ell,\beta_\ell)\in\B\) for
	\(1\leq\ell\leq k\). If \(\alpha_j=0\) for exactly one
	index \(j\), then
	\[
	u_1\odotop\cdots\odotop u_k
	=\left(0,\beta_j\prod_{\ell\ne j}\alpha_\ell\right)\in\B.
	\]
	If \(\alpha_j=0\) for at least two indices, then
	\(u_1\odotop\cdots\odotop u_k=(0,0)\).
\end{lemma}

\begin{proof}
	Repeated multiplication gives
	\[
	u_1\odotop\cdots\odotop u_k
	=\left(\prod_{\ell=1}^k\alpha_\ell,
	\sum_{j=1}^k\beta_j\prod_{\ell\ne j}\alpha_\ell\right).
	\]
	If exactly one \(\alpha_j\) is zero, only the \(j\)th term
	in the sum can be nonzero, and
	\[
	\left|\beta_j\prod_{\ell\ne j}\alpha_\ell\right|
	\leq\abs{\beta_j}\leq1.
	\]
	If at least two first coordinates are zero, every product in
	the displayed formula is zero.
\end{proof}

By Lemma~\ref{lem:zero-scalar}, every product containing a factor
with \(a_j=0\) belongs to \(\B\).
If \(\abs{a_j}=1\), the assumption \(s_j\in\B\) implies
\(b_j=0\), and
\[
z\odotop s_j=(a_jz_0,a_jz_1),\qquad
\norm{z\odotop s_j}_2=\norm z_2.
\]
Removing factors \(s_j=(a_j,0)\) with \(\abs{a_j}=1\)
from a product leaves its Euclidean norm unchanged.
These factors can therefore be omitted when checking
whether all products belong to \(\B\).
Thus it remains to consider products formed from elements
satisfying \(0<\abs{a_j}<1\).
If no such elements remain, every product either contains
a factor with zero first coordinate or has the form
\((\omega,0)\) with \(\abs\omega=1\).
Hence \(\mathcal R_S\subseteq\B\).
Otherwise, relabel the remaining elements as
\(S=\{s_1,\ldots,s_d\}\), with
\(0<\abs{a_j}<1\) for every \(j\). Put
\[
x_j=-\log\abs{a_j}>0,\quad c_j=\frac{b_j}{a_j},\quad
\bm x=(x_1,\ldots,x_d),\quad \bm c=(c_1,\ldots,c_d).
\]
For \(\bm n=(n_1,\ldots,n_d)\in\mathbb N_0^d\), define
\begin{equation}\label{eq:switching-Q}
	Q(\bm n)=e^{-2\bm n\cdot\bm x}
	\left(1+\abs{\bm n\cdot\bm c}^2\right).
\end{equation}
Set
\[
x_{\min}=\min_{1\leq j\leq d}x_j>0,\qquad
C_S=\max_{1\leq j\leq d}\abs{c_j}.
\]
Define
\begin{equation}\label{eq:switching-cutoff}
	N_*=\min\left\{N\in\mathbb N:
	N\geq x_{\min}^{-1},\
	e^{-2x_{\min}N}(1+C_S^2N^2)\leq1\right\}.
\end{equation}
The existence of \(N_*\) is established in the proof
of Proposition~\ref{prop:switching}.
The vector \(\bm n=\bm0\) gives the identity product and
\(Q(\bm0)=1\), so we only need to consider nonzero vectors.

For \(\bm n\in\mathbb N_0^d\), put
\(p(\bm n)=s_1^{\odotop n_1}\odotop\cdots\odotop
s_d^{\odotop n_d}\).
\begin{proposition}\label{prop:switching}
	All finite products of elements of \(S\) belong to \(\B\)
	if and only if
	\[
	Q(\bm n)\leq1
	\quad\text{for every }\bm n\in\mathbb N_0^d
	\text{ with }1\leq\norm{\bm n}_1<N_*.
	\]
\end{proposition}

\begin{proof}
	By associativity and commutativity, any finite product
	of elements of \(S\) in which \(s_j\) appears exactly
	\(n_j\) times for each \(j=1,\ldots,d\) equals
	\[
	s_1^{\odotop n_1}\odotop\cdots\odotop s_d^{\odotop n_d}
	=p(\bm n).
	\]
	The power formula gives
	\[
	p(\bm n)=
	\left(\prod_{j=1}^d a_j^{n_j},
	\left(\prod_{j=1}^d a_j^{n_j}\right)
	\sum_{j=1}^d n_jc_j\right).
	\]
	Consequently,
	\[
	\norm{p(\bm n)}_2^2
	=\left|\prod_{j=1}^d a_j^{n_j}\right|^2
	\left(1+\left|\sum_{j=1}^d n_jc_j\right|^2\right)
	=Q(\bm n).
	\]
	To bound \(Q(\bm n)\) in terms of the total number of
	factors, put \(k=\norm{\bm n}_1\). Then
	\[
	\begin{aligned}
		\bm n\cdot\bm x
		&=\sum_{j=1}^d n_jx_j
		\geq x_{\min}\sum_{j=1}^d n_j
		=x_{\min}k,\\
		\abs{\bm n\cdot\bm c}
		&\leq\sum_{j=1}^d n_j\abs{c_j}
		\leq C_S\sum_{j=1}^d n_j
		=C_Sk.
	\end{aligned}
	\]
	Therefore,
	\[
	Q(\bm n)
	\leq e^{-2x_{\min}k}(1+C_S^2k^2)
	=\Psi(k),
	\]
	where
	\[
	\Psi(t)=e^{-2x_{\min}t}(1+C_S^2t^2),
	\qquad t\geq0.
	\]
	Since \(x_{\min}>0\), both \(e^{-2x_{\min}t}\) and
	\(t^2e^{-2x_{\min}t}\) tend to zero as \(t\to\infty\).
	Thus \(\Psi(t)\to0\), which proves existence of \(N_*\).
	For \(t\geq x_{\min}^{-1}>0\),
	\[
	\begin{aligned}
		\frac{d}{dt}\log\Psi(t)
		&=-2x_{\min}
		+\frac2t\frac{C_S^2t^2}{1+C_S^2t^2}\\
		&\leq-2x_{\min}+\frac2t
		\leq0.
	\end{aligned}
	\]
	The inequalities follow from
	\(C_S^2t^2/(1+C_S^2t^2)\leq1\) and
	\(1/t\leq x_{\min}\).
	Since \(\Psi(t)>0\), we obtain \(\Psi'(t)\leq0\).
	Thus \(\Psi\) is nonincreasing on
	\([x_{\min}^{-1},\infty)\).
	Since \(N_*\geq x_{\min}^{-1}\) and
	\(\Psi(N_*)\leq1\), we obtain
	\[
	Q(\bm n)\leq\Psi(k)\leq\Psi(N_*)\leq1
	\qquad\text{whenever }k\geq N_*.
	\]
	If \(Q(\bm n)\leq1\) for every
	\(\bm n\in\mathbb N_0^d\) satisfying
	\(1\leq\norm{\bm n}_1<N_*\), then
	\(Q(\bm n)\leq1\) for every nonzero \(\bm n\).
	Together with \(Q(\bm0)=1\), we can obtain
	\(p(\bm n)\in\B\) for all \(\bm n\).
	
	Conversely, if \(p(\bm n)\in\B\) for every
	\(\bm n\in\mathbb N_0^d\), then
	\[
	Q(\bm n)=\norm{p(\bm n)}_2^2\leq1
	\qquad\text{whenever }1\leq\norm{\bm n}_1<N_*.
	\]
\end{proof}

The number of inequalities \(Q(\bm n)\leq1\)
that must be checked in the finite criterion is
\[
L_d(N_*)=
\#\{\bm n\in\mathbb N_0^d:1\leq\norm{\bm n}_1<N_*\}
=\binom{N_*+d-1}{d}-1.
\]
For fixed \(d\), the number of inequalities to be checked
grows polynomially in \(N_*\), with leading term
\(N_*^d/d!\). More precisely,
\[
\lim_{N_*\to\infty}
\frac{L_d(N_*)}{N_*^d/d!}=1.
\]
Commutativity allows products with the same multiplicities to
be considered together. The number \(L_d(N_*)\) nevertheless
grows with \(d\) and \(N_*\). 

Consider \(S=\{p,q\}\), where \(p,q\) are the elements
in \eqref{eq:pq}.
All powers of either element belong to \(\B\).
The finite criterion nevertheless detects a mixed
product outside \(\B\), as summarized in
Table~\ref{tab:finite-product-example}.

\begin{table}[H]
	\caption{The finite product criterion for the pair \(p,q\).}
	\label{tab:finite-product-example}
	\centering
	\begin{tabular}{@{}ll@{}}
		\toprule
		Quantity & Value\\
		\midrule
		\(x_{\min}\) & \(\log(4/3)\approx0.2876820725\)\\
		\(C_S\) & \(\sqrt2\)\\
		\(N_*\) & \(9\)\\
		Number of inequalities \(L_2(9)\) & \(44\)\\
		Maximizing index \((n_1,n_2)\) & \((1,1)\)\\
		Maximum of \(Q(\bm n)\) & \(1.0256060601\)\\
		\bottomrule
	\end{tabular}
\end{table}

Here
\[
Q(n_1,n_2)
=\left(\frac9{16}\right)^{n_1}3^{-n_2}
\left[1+\left(\frac7{10}n_1+\sqrt2\,n_2\right)^2\right],
\]
and the maximum is taken over
\(1\leq n_1+n_2<9\).
In particular,
\[
Q(1,1)=\frac{1047}{1600}+\frac{21\sqrt2}{80}>1.
\]
Thus the criterion fails because of the mixed product
\(p\odotop q\), although every pure power remains in \(\B\).

For this example,
\[
\frac14<\log(4/3)<\frac13,
\qquad
\lceil x_{\min}^{-1}\rceil=4,
\qquad
\Psi(N)=\left(\frac9{16}\right)^N(1+2N^2).
\]
Since \(\Psi\) is nonincreasing for \(N\geq4\), the
exact inequalities
\[
\Psi(8)=\frac{5553027009}{4294967296}>1,
\qquad
\Psi(9)=\frac{63149539707}{68719476736}<1
\]
give \(N_*=9\).
Exact comparison in \(\mathbb Q(\sqrt2)\) over
\(1\leq n_1+n_2<9\) shows that \((1,1)\) is the
only index with \(Q(n_1,n_2)>1\).
Equality holds exactly at \((0,1)\) and \((0,2)\).
The displayed decimals are rounded for presentation.
All products with at least nine factors belong to
\(\B\) by the analytic tail estimate.

We now consider prescribed dual numbers
\(s_j(\theta)=(a_j(\theta),b_j(\theta))\)
depending on a real parameter vector \(\theta\).
Under the assumptions below, a single integer \(N\)
can be chosen so that every product with at least \(N\)
factors belongs to \(\B\) throughout a neighborhood
of a fixed parameter value.
The condition for all finite products to belong to
\(\B\) then reduces to inequalities with the same
finite index set throughout this neighborhood.

Let \(\Theta\subseteq\mathbb R^m\) be open, and let
\(a_j,b_j:\Theta\to\C\) be continuous for
\(1\leq j\leq d\). Fix \(\theta_0\in\Theta\) such that
\(0<\abs{a_j(\theta_0)}<1\) for every \(j\).
Since each \(a_j\) is continuous and
\(0<\abs{a_j(\theta_0)}<1\) for \(1\leq j\leq d\),
there exist an open neighborhood \(U\) of \(\theta_0\),
with compact closure \(\overline U\subseteq\Theta\), and
constants \(a_-,a_+\) such that
\[
0<a_-\leq\abs{a_j(\theta)}
\leq a_+<1
\qquad
(\theta\in\overline U,\ 1\leq j\leq d).
\]
Set
\[
\begin{aligned}
	x_j(\theta)&=-\log\abs{a_j(\theta)},&
	c_j(\theta)&=\frac{b_j(\theta)}{a_j(\theta)},\\
	\alpha&=\min_{\substack{\theta\in\overline U\\1\leq j\leq d}}
	x_j(\theta)>0,&
	C&=\max_{\substack{\theta\in\overline U\\1\leq j\leq d}}
	\abs{c_j(\theta)}<\infty.
\end{aligned}
\]
Since \(a_j\) does not vanish on \(\overline U\),
the function \(c_j=b_j/a_j\) is continuous there.
Compactness of \(\overline U\) therefore implies that
the maximum defining \(C\) is finite. Moreover,
\[
x_j(\theta)\geq-\log a_+>0
\qquad
(\theta\in\overline U,\ 1\leq j\leq d),
\]
so \(\alpha>0\).

Since \(e^{-2\alpha t}(1+C^2t^2)\to0\), choose a fixed
\(N\in\mathbb N\) satisfying
\begin{equation}\label{eq:uniform-cutoff}
	N\geq\alpha^{-1},\qquad e^{-2\alpha N}(1+C^2N^2)\leq1.
\end{equation}
For \(\bm n\in\mathbb N_0^d\), define
\begin{equation}\label{eq:parameter-Q}
	Q_\theta(\bm n)=e^{-2\sum_jn_jx_j(\theta)}
	\left(1+\left|\sum_jn_jc_j(\theta)\right|^2\right).
\end{equation}
Put
\[
I_N=\{\bm n\in\mathbb N_0^d:1\leq\norm{\bm n}_1<N\},
\qquad
\Phi(\theta)=\max_{\bm n\in I_N}Q_\theta(\bm n),
\]
with \(\Phi\equiv0\) when \(I_N=\varnothing\).

The uniform cutoff \(N\) gives the following finite
criterion throughout \(U\).

\begin{corollary}\label{cor:uniform-criterion}
	For each \(\theta\in U\), all finite products of
	\[
	(a_1(\theta),b_1(\theta)),\ldots,
	(a_d(\theta),b_d(\theta)),
	\]
	including the identity, belong to \(\B\) if and only if
	\[
	\Phi(\theta)\leq1.
	\]
\end{corollary}

\begin{proof}
	Fix \(\theta\in U\). For \(\bm n\in\mathbb N_0^d\),
	put \(k=\norm{\bm n}_1\).
	Since \(x_j(\theta)\geq\alpha\) and
	\(\abs{c_j(\theta)}\leq C\), we have
	\[
	\begin{aligned}
		\sum_{j=1}^d n_jx_j(\theta)
		&\geq\alpha\sum_{j=1}^d n_j=\alpha k,\\
		\left|\sum_{j=1}^d n_jc_j(\theta)\right|
		&\leq\sum_{j=1}^d n_j\abs{c_j(\theta)}
		\leq Ck.
	\end{aligned}
	\]
	Consequently,
	\[
	Q_\theta(\bm n)\leq e^{-2\alpha k}(1+C^2k^2).
	\]
	Let
	\[
	\Psi(t)=e^{-2\alpha t}(1+C^2t^2).
	\]
	For \(t\geq\alpha^{-1}\),
	\[
	\frac{d}{dt}\log\Psi(t)
	=-2\alpha+\frac{2C^2t}{1+C^2t^2}
	\leq-2\alpha+\frac2t\leq0.
	\]
	Since \(\Psi(t)>0\), it follows that \(\Psi\) is
	nonincreasing on \([\alpha^{-1},\infty)\).
	By \eqref{eq:uniform-cutoff},
	\(N\geq\alpha^{-1}\) and \(\Psi(N)\leq1\).
	
	Thus, if \(k\geq N\), then
	\[
	Q_\theta(\bm n)
	\leq\Psi(k)\leq\Psi(N)\leq1.
	\]
	Since \(Q_\theta(\bm0)=1\), all finite products belong
	to \(\B\) if and only if
	\[
	Q_\theta(\bm n)\leq1
	\qquad\text{for every }\bm n\in I_N.
	\]
	If \(I_N\ne\varnothing\), these inequalities are
	equivalent to \(\Phi(\theta)\leq1\).
	If \(I_N=\varnothing\), the estimate for \(k\geq N\)
	already ensures that every nonempty product belongs
	to \(\B\), and \(\Phi(\theta)=0\leq1\) by definition.
\end{proof}

\subsection{An invariant error region for an insulin-delay subsystem}
\label{subsec:insulin-delay}

In this subsection, we use \(M_E^{\max}\cap\mathbb R^2\)
as an invariant error region for an insulin-delay subsystem.
M. Nagumo \cite{Nagumo1942} established a classical tangency
criterion for positive invariance.
Invariant sets under bounded disturbances and their analysis
by Minkowski functionals are treated in
\cite{Blanchini1999,BlanchiniMiani2015}.
For the error region considered here, we compute its
Minkowski functional explicitly and determine the exact
input bound for robust positive invariance.

The model equations are taken from
X.~Chen, S.~Dutta and S.~Sankaranarayanan
\cite[Section~3.2, p.~82]{Chen2017} and are given by
\begin{equation}\label{eq:insulin-model}
	\dot I_1=-r(I_1-\gamma I_p),\quad
	\dot I_d=-r(I_d-I_1),\quad
	r=0.0046\,\mathrm{min}^{-1},\quad \gamma=18.2129.
\end{equation}
Here \(I_p\) denotes the plasma insulin concentration.
It is treated as an external input to this subsystem.
The variable \(I_1\) is the insulin concentration in the
first compartment, and \(I_d\) represents the delayed
response to \(I_1\).
For prescribed initial values and inputs \(I_p\) and
\(\bar I_p\), let \((I_1,I_d)\) and
\((\bar I_1,\bar I_d)\) be the corresponding solutions
of \eqref{eq:insulin-model}, respectively.
To express the errors in dimensionless form, choose a
common scale \(L>0\) with the same units as \(I_1\)
and \(I_d\), and define
\[
e_0=\frac{I_1-\bar I_1}{L},\qquad
e_1=\frac{I_d-\bar I_d}{L},\qquad
v=\frac{\gamma(I_p-\bar I_p)}{L}.
\]
\begin{lemma}\label{lem:insulin-error-flow}
	The normalized error \(e=(e_0,e_1)\) satisfies
	\begin{equation}\label{eq:insulin-error}
		\dot e_0=r(-e_0+v),\qquad
		\dot e_1=r(e_0-e_1),
	\end{equation}
	or, equivalently,
	\[
	\dot e=(-r,r)\odotop e+(rv,0).
	\]
	For \(v\equiv0\), the solution with initial value
	\(e(0)\in\mathbb R^2\) is given by
	\[
	e(t)=s(t)\odotop e(0),
	\qquad s(t)=e^{-rt}(1,rt),
	\qquad t\geq0.
	\]
	Moreover,
	\[
	s(t+u)=s(t)\odotop s(u)
	\qquad(t,u\geq0).
	\]
\end{lemma}

\begin{proof}
	Subtracting the equations for \((\bar I_1,\bar I_d)\)
	from the corresponding equations for \((I_1,I_d)\) gives
	\[
	\begin{aligned}
		\dot I_1-\dot{\bar I}_1
		&=-r(I_1-\bar I_1)+r\gamma(I_p-\bar I_p),\\
		\dot I_d-\dot{\bar I}_d
		&=r(I_1-\bar I_1)-r(I_d-\bar I_d).
	\end{aligned}
	\]
	Since \(L\) is constant, division by \(L\) and the
	definitions of \(e_0,e_1,v\) yield \eqref{eq:insulin-error}.
	
	For \(e=(e_0,e_1)\), the definition of dual-number
	multiplication implies
	\[
	\begin{aligned}
		(-r,r)\odotop e+(rv,0)
		&=(-re_0,re_0-re_1)+(rv,0)\\
		&=\bigl(r(-e_0+v),r(e_0-e_1)\bigr).
	\end{aligned}
	\]
	Since \(\dot e=(\dot e_0,\dot e_1)\), the two equations
	in \eqref{eq:insulin-error} can therefore be written as
	\[
	\dot e=(-r,r)\odotop e+(rv,0).
	\]
	For \(v\equiv0\), the error equations \eqref{eq:insulin-error} become
	\begin{equation}\label{1}
		\dot e_0=-re_0,\qquad
		\dot e_1+re_1=re_0.
	\end{equation}
	Multiplying the first equation in \eqref{1} by \(e^{rt}\) gives
	\[
	\frac{d}{dt}\bigl(e^{rt}e_0(t)\bigr)=0.
	\]
	Hence \(e^{rt}e_0(t)=e_0(0)\), and therefore
	\[
	e_0(t)=e^{-rt}e_0(0).
	\]
	Substitution into the second equation in \eqref{1} yields
	\[
	\dot e_1+re_1=re^{-rt}e_0(0),
	\qquad
	\frac{d}{dt}\bigl(e^{rt}e_1(t)\bigr)=re_0(0).
	\]
	Integrating from \(0\) to \(t\), we obtain
	\[
	e^{rt}e_1(t)-e_1(0)=rt\,e_0(0),
	\]
	so that
	\[
	e_1(t)=e^{-rt}\bigl(e_1(0)+rt\,e_0(0)\bigr).
	\]
	
	Define
	\[
	s(t)=e^{-rt}(1,rt),\qquad t\geq0.
	\]
	Here \(e(0)=(e_0(0),e_1(0))\) is the initial error vector. By the definition of dual-number multiplication,
	\[
	\begin{aligned}
		s(t)\odotop e(0)
		&=(e^{-rt},rt\,e^{-rt})
		\odotop(e_0(0),e_1(0))\\
		&=\bigl(e^{-rt}e_0(0),
		e^{-rt}(e_1(0)+rt\,e_0(0))\bigr)\\
		&=e(t).
	\end{aligned}
	\]
	Moreover, for \(t,u\geq0\),
	\[
	\begin{aligned}
		s(t)\odotop s(u)
		&=e^{-r(t+u)}(1,rt)\odotop(1,ru)\\
		&=e^{-r(t+u)}(1,r(t+u))\\
		&=s(t+u).
	\end{aligned}
	\]
\end{proof}
Let \(K=M_E^{\max}\cap\mathbb R^2\). It follows from Theorem~\ref{thm:explicit-maximal} and
\eqref{eq:common-inner-ball} that \(K\) is a convex body
symmetric about the origin, with \(0\) in its interior. Define the Minkowski functional of \(K\) by
\[
p_K(e)=\inf\{q>0:e\in qK\},
\qquad e\in\mathbb R^2,
\]
where \(qK=\{qx:x\in K\}\).
Thus \(p_K\) is a norm on \(\mathbb R^2\)
whose closed unit ball is \(K\).

\begin{lemma}\label{lem:insulin-norm}
	The Minkowski functional \(p_K\) is given by
	\begin{equation}\label{eq:insulin-gauge}
		p_K(e_0,e_1)=
		\begin{cases}
			\displaystyle
			\sqrt{e_0^2+\frac{2\sqrt2}{3}\abs{e_0e_1}
				+\frac13e_1^2},
			&\abs{e_1}\leq\sqrt2\abs{e_0},\\[2mm]
			\sqrt{e_0^2+e_1^2},
			&\abs{e_1}\geq\sqrt2\abs{e_0}.
		\end{cases}
	\end{equation}
	It satisfies
	\begin{equation}\label{eq:insulin-submultiplicative}
		p_K(a\odotop b)\leq p_K(a)p_K(b)
		\qquad(a,b\in\mathbb R^2),
	\end{equation}
	and \(p_K(v,0)=\abs v\) for every \(v\in\mathbb R\).
\end{lemma}

\begin{proof}
	For \(e\ne0\), put \(q=p_K(e)>0\).
	Then \(e/q\in\partial K\), so
	\eqref{eq:explicit-maximal} yields
	\[
	\frac{\abs{e_1}}{q}
	=\rho_{\max}\!\left(\frac{\abs{e_0}}{q}\right).
	\]
	If \(\abs{e_0}/q\leq1/\sqrt3\), the first expression
	in \eqref{eq:rho-max} implies that
	\[
	q^2=e_0^2+e_1^2,
	\]
	and this case corresponds to
	\(\abs{e_1}\geq\sqrt2\abs{e_0}\).
	If \(\abs{e_0}/q\geq1/\sqrt3\), the second expression shows that
	\[
	\abs{e_1}+\sqrt2\abs{e_0}
	=\sqrt{3q^2-e_0^2}.
	\]
	Squaring and rearranging yields
	\[
	q^2=e_0^2+\frac{2\sqrt2}{3}\abs{e_0e_1}
	+\frac13e_1^2,
	\]
	with the corresponding condition
	\(\abs{e_1}\leq\sqrt2\abs{e_0}\).
	Together with \(p_K(0,0)=0\), these calculations give
	\eqref{eq:insulin-gauge}.
	
	In particular, \(p_K(v,0)=\abs v\) for \(v\in\mathbb R\).
	Since \(K=M_E^{\max}\cap\mathbb R^2\), its Minkowski
	functional satisfies
	\(p_K=p_{M_E^{\max}}|_{\mathbb R^2}\).
	Thus \eqref{eq:insulin-submultiplicative} follows from
	Corollary~\ref{cor:maximal-algebra-norm}.
\end{proof}
Let \(W\geq0\). The set \(K\) is called
\emph{robust positively invariant} under the input bound
\(W\) if, for every initial value \(e(0)\in K\) and every
real measurable input \(v\) satisfying
\(\abs{v(t)}\leq W\) for almost every \(t\geq0\),
the corresponding solution of \eqref{eq:insulin-error}
satisfies
\[
e(t)\in K\qquad\text{for all }t\geq0.
\]
This is the standard notion of robust positive invariance
\cite{Blanchini1999}. The following theorem determines the
largest value of \(W\) for which \(K\) has this property.

\begin{theorem}\label{thm:insulin-invariance}
	Let \(r>0\), \(W\geq0\), and \(\kappa=1-\sqrt2/3\).
	For the system
	\[
	\dot e=(-r,r)\odotop e+(rv,0),
	\]
	\(K\) is robust positively invariant under the input
	bound \(W\) if and only if \(W\leq\kappa\).
	
	For every initial value \(e(0)\in\mathbb R^2\) and
	every real measurable input \(v\) satisfying
	\(\abs{v(t)}\leq W\) for almost every \(t\geq0\),
	the corresponding solution \(e(t)\) of the error
	system satisfies
	\begin{equation}\label{eq:insulin-estimate}
		p_K(e(t))
		\leq e^{-r\kappa t}p_K(e(0))
		+\frac{W}{\kappa}\bigl(1-e^{-r\kappa t}\bigr),
		\qquad t\geq0.
	\end{equation}
\end{theorem}

\begin{proof}
	For \(0\leq h\leq\sqrt2/r\), the coordinates of
	\[
	s(h)=e^{-rh}(1,rh)
	\]
	satisfy
	\[
	rh\,e^{-rh}\leq\sqrt2\,e^{-rh}.
	\]
	Substituting \(s(h)=(e^{-rh},rh\,e^{-rh})\) into the
	first expression in \eqref{eq:insulin-gauge} gives
	\[
	\begin{aligned}
		p_K(s(h))
		&=\sqrt{e^{-2rh}
			+\frac{2\sqrt2}{3}rh\,e^{-2rh}
			+\frac13r^2h^2e^{-2rh}}\\
		&=e^{-rh}
		\sqrt{1+\frac{2\sqrt2}{3}rh+\frac13r^2h^2}.
	\end{aligned}
	\]
	The Taylor expansions at \(h=0\) give
	\[
	\begin{aligned}
		e^{-rh}
		&=1-rh+O(h^2),\\
		\sqrt{1+\frac{2\sqrt2}{3}rh+\frac13r^2h^2}
		&=1+\frac{\sqrt2}{3}rh+O(h^2).
	\end{aligned}
	\]
	Consequently,
	\[
	\begin{aligned}
		p_K(s(h))
		&=(1-rh+O(h^2))
		\left(1+\frac{\sqrt2}{3}rh+O(h^2)\right)\\
		&=1-r\left(1-\frac{\sqrt2}{3}\right)h+O(h^2)\\
		&=1-r\kappa h+O(h^2),
		\qquad h\to0^+.
	\end{aligned}
	\]
	
	Fix \(t>0\). For all sufficiently large \(n\),
	\(t/n\leq\sqrt2/r\). The semigroup identity
	\[
	s(t)=s(t/n)^{\odotop n}
	\]
	and \eqref{eq:insulin-submultiplicative} imply
	\[
	\begin{aligned}
		p_K(s(t))
		&\leq p_K(s(t/n))^n \\
		&=
		\left(
		1-\frac{r\kappa t}{n}+O(n^{-2})
		\right)^n.
	\end{aligned}
	\]
	Letting \(n\to\infty\), we obtain
	\[
	p_K(s(t))\leq e^{-r\kappa t}.
	\]
	The same inequality holds at \(t=0\), since
	\(p_K(s(0))=p_K(1,0)=1\).

	Since the homogeneous system has constant coefficients,
	its solution with value \(z\) at time \(u\) is
	\(s(t-u)\odotop z\) for \(t\geq u\).
	The variation-of-constants formula for
	\eqref{eq:insulin-error} therefore takes the form
	\[
	e(t)=s(t)\odotop e(0)
	+r\int_0^t s(t-u)\odotop(v(u),0)\,du.
	\]
	
	To verify this formula, denote its right-hand side by
	\(F(t)\). Since \(v\) is measurable and bounded almost
	everywhere, and \(s\) is continuously differentiable,
	\(F\) is locally absolutely continuous.
	Moreover, \(s(0)=(1,0)\), so
	\[
	F(0)=s(0)\odotop e(0)=e(0).
	\]
	Direct differentiation of \(s(t)=e^{-rt}(1,rt)\) gives
	\[
	\begin{aligned}
		s'(t)
		&=\bigl(-re^{-rt},r(1-rt)e^{-rt}\bigr)\\
		&=(-r,r)\odotop s(t).
	\end{aligned}
	\]
	
	In the integral defining \(F(t)\), both the upper limit
	and the factor \(s(t-u)\) depend on \(t\).
	Differentiating almost everywhere yields
	\[
	\begin{aligned}
		F'(t)
		&=s'(t)\odotop e(0)
		+r\,s(0)\odotop(v(t),0)\\
		&\quad
		+r\int_0^t s'(t-u)\odotop(v(u),0)\,du.
	\end{aligned}
	\]
	Here the term \(r\,s(0)\odotop(v(t),0)\) arises from
	the upper limit of integration.
	Substituting
	\(s'(\tau)=(-r,r)\odotop s(\tau)\) and \(s(0)=(1,0)\)
	gives
	\[
	F'(t)=(-r,r)\odotop F(t)+(rv(t),0)
	\quad\text{almost everywhere}.
	\]
	Thus \(F\) satisfies \eqref{eq:insulin-error} with the
	initial value \(F(0)=e(0)\). Uniqueness of the solution
	implies \(F(t)=e(t)\) for all \(t\geq0\).
	
	Using the triangle inequality, the integral inequality for
	the norm \(p_K\), \eqref{eq:insulin-submultiplicative}, and
	\(p_K(v(u),0)=\abs{v(u)}\), we obtain
	\begin{align*}
		p_K(e(t))
		&\leq
		p_K(s(t)\odotop e(0))
		+r\int_0^t
		p_K\bigl(s(t-u)\odotop(v(u),0)\bigr)\,du\\
		&\leq
		p_K(s(t))p_K(e(0))
		+r\int_0^t
		p_K(s(t-u))\abs{v(u)}\,du\\
		&\leq
		e^{-r\kappa t}p_K(e(0))
		+rW\int_0^t e^{-r\kappa(t-u)}\,du\\
		&=
		e^{-r\kappa t}p_K(e(0))
		+\frac{W}{\kappa}
		\bigl(1-e^{-r\kappa t}\bigr).
	\end{align*}
	This proves \eqref{eq:insulin-estimate}.
	
	If \(e(0)\in K\), then \(p_K(e(0))\leq1\).
	When \(W\leq\kappa\), \eqref{eq:insulin-estimate} implies
	\[
	\begin{aligned}
		p_K(e(t))
		&\leq
		e^{-r\kappa t}
		+\frac{W}{\kappa}
		\bigl(1-e^{-r\kappa t}\bigr)\\
		&\leq
		e^{-r\kappa t}
		+1-e^{-r\kappa t}
		=1.
	\end{aligned}
	\]
	Hence \(e(t)\in K\) for every \(t\geq0\).
	
	To prove necessity, take the boundary point
	\(e(0)=(1,0)\in\partial K\) and the admissible
	constant input \(v(t)\equiv W\). By equation~\eqref{eq:insulin-error}, we have 
	\begin{equation}\label{eq:insulin-initial-expansion}
		e_0(t)=1+r(W-1)t+O(t^2),
		\qquad
		e_1(t)=rt+O(t^2).
	\end{equation}
	For all sufficiently small \(t>0\), these expansions imply
	\[
	e_0(t)>0,\qquad
	e_1(t)>0,\qquad
	e_1(t)\leq\sqrt2 e_0(t).
	\]
	The first expression in \eqref{eq:insulin-gauge} therefore shows that
	\[
	p_K(e(t))
	=1+r(W-\kappa)t+O(t^2).
	\]
	If \(W>\kappa\), then \(p_K(e(t))>1\) for all sufficiently
	small \(t>0\). Thus the solution leaves \(K\), and \(K\)
	is not robust positively invariant.
\end{proof}
The threshold \(\kappa\) is exact for robust positive
invariance of \(K\). The estimate
\eqref{eq:insulin-estimate} is an upper bound on the
solution norm and need not be attained.

For comparison, we consider the Toeplitz set \(\Tset\)
introduced earlier in this paper, restricted to real
coordinates. The resulting region and its Minkowski
functional are
\[
K_T=\Tset\cap\mathbb R^2
=\{(e_0,e_1):e_0^2+\abs{e_1}\leq1\},\qquad
p_{K_T}(e_0,e_1)=\frac{\sqrt{4e_0^2+e_1^2}+\abs{e_1}}2.
\]
The same semigroup argument, applied to
\(p_{K_T}(s(h))=1-rh/2+O(h^2)\), gives
\[
p_{K_T}(e(t))
\leq e^{-rt/2}p_{K_T}(e(0))
+2W(1-e^{-rt/2}),
\quad t\geq0.
\]
Thus \(K_T\) is robust positively invariant when
\(W\leq1/2\).
For the initial value \((1,0)\) and the constant
input \(v(t)=W\),
\eqref{eq:insulin-initial-expansion} gives
\[
p_{K_T}(e(t))=1+r(W-1/2)t+O(t^2).
\]
Hence \(K_T\) is robust positively invariant under
all measurable inputs satisfying
\(\abs{v(t)}\leq W\) almost everywhere
if and only if \(W\leq1/2\).
Thus the passage from \(K_T\) to the larger region \(K\)
increases the exact admissible input bound from
$1/2$ to $1-\sqrt{2}/3$.

\begin{figure}[H]
	\centering
	\includegraphics[width=0.8\linewidth]{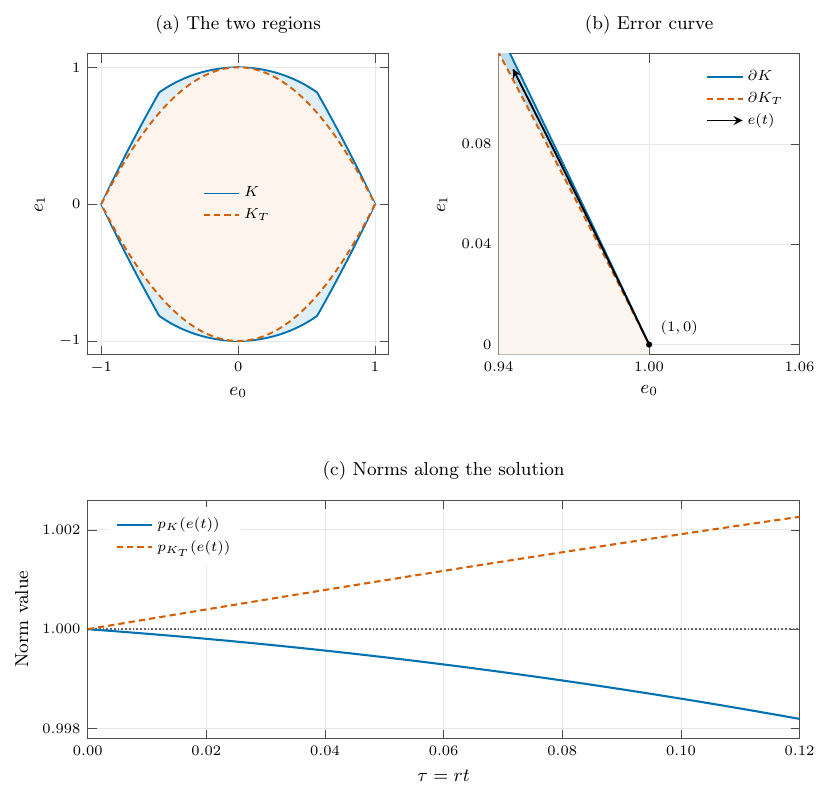}
	\caption{The regions \(K_T\subsetneq K\) and the solution
		with \(e(0)=(1,0)\) and \(v(t)\equiv0.52\).
		(a) The two regions, with \(K\) bounded by the solid
		blue curve and \(K_T\) by the dashed orange curve.
		(b) The error curve near \((1,0)\), for
		\(0\leq\tau=rt\leq0.12\).
		The blue shaded strip lies in \(K\setminus K_T\),
		and the arrow indicates increasing time.
		(c) The norms \(p_K(e(t))\) and \(p_{K_T}(e(t))\)
		over the same time interval.
		The horizontal dotted line is the boundary level \(1\).}
	\label{fig:insulin-regions}
\end{figure}

To illustrate the different admissible input bounds for
\(K_T\) and \(K\), take \(e(0)=(1,0)\) and the constant
input \(v(t)\equiv W=0.52\), so that
\(1/2<W<\kappa\).
The solution of \eqref{eq:insulin-error} is
\[
\begin{aligned}
	e_0(t)&=W+(1-W)e^{-rt},\\
	e_1(t)&=W(1-e^{-rt})+(1-W)rt\,e^{-rt}.
\end{aligned}
\]
By Theorem~\ref{thm:insulin-invariance},
\(e(t)\in K\) for every \(t\geq0\).
Since
\[
p_{K_T}(e(t))=1+r(W-1/2)t+O(t^2),
\]
we have \(e(t)\notin K_T\) for all sufficiently small
\(t>0\).
Figure~\ref{fig:insulin-regions} shows the regions \(K_T\) and
\(K\), the error curve near \((1,0)\), and the two norm
values for \(0\leq\tau=rt\leq0.12\).
A norm value at most \(1\) indicates that the solution
lies in the corresponding region.

In the original variables, the sufficient input condition is
\begin{equation}\label{eq:insulin-physical-bound}
	\abs{I_p(t)-\bar I_p(t)}
	\leq\frac{L\kappa}{\gamma}
	\quad\text{for almost every }t\geq0.
\end{equation}
Here
\[
\frac{\kappa}{\gamma}
=\frac{1-\sqrt2/3}{18.2129}
\approx0.0290231363.
\]
The scale \(L\) specifies the size of the error region \(LK\).
The bound in \eqref{eq:insulin-physical-bound} is therefore
relative to this chosen region. No particular physical
value of \(L\) is selected here.

Under this condition, if
\[
\bigl(I_1(0)-\bar I_1(0),
I_d(0)-\bar I_d(0)\bigr)\in LK,
\]
then
\[
\bigl(I_1(t)-\bar I_1(t),
I_d(t)-\bar I_d(t)\bigr)\in LK,
\quad t\geq0.
\]
Applying the conclusion to a complete closed-loop
model requires a separate verification of
\eqref{eq:insulin-physical-bound} for that model.

For the same error system, the real Euclidean unit ball
\(B_{\mathbb R}=\B\cap\mathbb R^2\) admits a larger
input bound. At a boundary point \((a,b)\), the outward
normal is \((a,b)\). For \(\abs v\leq W\), the vector
field in \eqref{eq:insulin-error} satisfies
\[
\frac{a(-ra+rv)+b(ra-rb)}{r}
=-1+ab+av
\leq-1+\abs a\abs b+W\abs a.
\]
The boundary tangency condition
\cite{Nagumo1942,BlanchiniMiani2015} gives the exact
threshold
\[
W_B=\min_{0<a\leq1}
\left(\frac1a-\sqrt{1-a^2}\right).
\]
The expression tends to infinity as \(a\downarrow0\),
and its derivative is
\[
-\frac1{a^2}+\frac{a}{\sqrt{1-a^2}},
\qquad 0<a<1.
\]
This derivative changes sign from negative to positive
at \(a_*=\sqrt y\), where \(y\in(0,1)\) is the unique
root of \(y^3+y-1=0\). Thus
\[
W_B=\frac1{\sqrt y}-\sqrt{1-y}
\approx0.6469836322.
\]
If \(W>W_B\), take
\(e(0)=(\sqrt y,\sqrt{1-y})\) in the first quadrant
and \(v(t)\equiv W\). The outward normal component
at \(t=0\) is \(r\sqrt y\,(W-W_B)>0\), so the
solution leaves \(B_{\mathbb R}\).
Consequently,
\[
\frac12<\kappa<W_B.
\]
The improvement from \(K_T\) to \(K\) concerns the two
specified multiplicatively closed regions. It does not
optimize the input bound over all invariant regions
contained in the Euclidean unit ball.

\section*{Author contributions}
Cailing Yao: Conceptualization, Writing-Original Draft, Formal Analysis, Methodology; Bingzhe Hou: Conceptualization, Formal Analysis, Methodology; Chao Xie: Methodology, Writing-review $\&$ editing; Yuanjun Zou: Methodology, Formal Analysis, Writing-review $\&$ editing.

\section*{Use of Generative-AI tools declaration}
The authors declare that Artiffcial Intelligence (AI) tools were used solely for language checking of this manuscript. No AI tools were used to generate research content, data analysis, or original text.

\section*{Acknowledgments}
We thank the editor and the reviewers for their dedicated time and constructive comments. 

\section*{Conflict of interest}
The corresponding authors declare that there are no conflicts of interest in this work.

\end{document}